\documentclass[12pt]{amsart}
\usepackage{amssymb, latexsym, amsthm} 

\newtheorem{thm}{Theorem}[section] 
\newtheorem{defn}[thm]{Definition}

\newtheorem{cor}[thm]{Corollary}

\newtheorem{exmpl}[thm]{Example}
\newtheorem{prop}[thm]{Proposition}
\newtheorem{rem}[thm]{Remark}

\def\ie{{i.e.}} 
\def\eg{{e.g.}}

\def\s{\sigma}
\def\ta{\tau}

\def\C{{\mathbb C}}

\def\F{{\mathbb F}}

\def\P{{\mathbb{P}}}

\def\T{{\mathbb T}}

\def\Z{{\mathbb Z}}

\def\ra{{\rightarrow}}

\def\llra{{\longrightarrow}}

\def\sub{\subseteq}
\def\({\left(}
\def\){\right)}
\def\isom{{\cong}}

\def\suchthat{{\,:\,}}

\def\maptoSn{{\phi}}

\DeclareMathOperator{\ab}{ab}

\DeclareMathOperator{\Ker}{Ker}
\renewcommand{\Im}{\operatorname{Im}}

\newcommand\til[1]{{\widetilde{#1}}}
\newcommand\tilmap[1][u]{{{#1} \mapsto \til{#1}}}
\newcommand\Cox[1]{{\operatorname{C}(#1)}}
\newcommand\CoxY[1]{{\operatorname{C}_{\operatorname{Y}}(#1)}}

\newcommand\eq[1]{{(\ref{#1})}}
\newcommand\eqs[3][--]{{\eq{#2}#1\eq{#3}}}
\newcommand\Eq[1]{{Equation \eq{#1}}}

\newcommand\Cref[1]{{Corollary \ref{#1}}}
\newcommand{\set}[1]{{\{#1\}}}
\newcommand{\power}[1]{{\left|{#1}\right|}}
\newcommand\defin[1]{{\it{#1}}}

\newcommand\ideal[1]{{\left<{#1}\right>}}
\newcommand\sg[1]{{\ideal{#1}}}

\long\def\forget#1\forgotten{}
\newcommand\explain[2]{{\stackrel{{{#2}}}{#1}}}
\newcommand\eqY{{\stackrel{Y}{=}}}

\newcommand\restricted[1]{{|_{#1}}} 
\newcommand\semidirect[2]{{{#1}\ltimes {#2}}} 

\newcommand\Astar[1][t,n]{{A^{*}_{#1}}}
\newcommand\Fstar[1][t,n]{{F^{*}_{#1}}}
\newcommand\Tsar{{\operatorname{Ts}^{*}(\Gamma)}}
\def\co{{\,{:}\,}}

\newif \iffurther 
 \furtherfalse

\newif\ifXY 
\XYtrue     

\ifXY
\usepackage{xy}
\fi \ifXY \xyoption{all} \fi

\begin{document}

\title{Coxeter Covers of the Symmetric Groups}

\author{Louis Rowen, Mina Teicher and Uzi Vishne}
\address{
L.H.~Rowen, M.~Teicher, Department of Mathematics, Bar-Ilan
University, Ramat-Gan 52900, Israel } \email{rowen@math.biu.ac.il,
teicher@math.biu.ac.il}

\address
{ Uzi Vishne, Department of Mathematics, Yale University, 10
Hillhouse Ave. New-Haven CT 06520, USA } \email{vishne@math.biu.ac.il}

\date{15 Nov. 2003} 

\thanks{The third named author was partially supported by the Fulbright Visiting Scholar
Program, United States Department of State}

\begin{abstract}
We study Coxeter groups from which there is a natural map onto a
symmetric group. Such groups have natural quotient groups related
to presentations of the symmetric group on an arbitrary set $T$ of
transpositions. These quotients, denoted here by $\CoxY{T}$, are a
special type of the generalized Coxeter groups defined in
\cite{CST}, and also arise in the computation of certain
invariants of surfaces.

We use a surprising action of $S_n$ on the kernel of the
surjection $\CoxY{T} \ra S_n$ to show that this kernel embeds in
the direct product of $n$ copies of the free group $\pi_1(T)$
(with the exception of $T$ being the full set of transpositions in
$S_4$). As a result, we show that the groups $\CoxY{T}$ are either
virtually Abelian or contain a non-Abelian free subgroup.
\end{abstract}

\maketitle

\section{Introduction}\label{intro}

The symmetric group on $n$ letters is generated by the transpositions
$s_i = (i\,i+1)$, $i = 1,\dots,n-1$. These generators satisfy the well known relations
$s_i^2 = 1$, $s_is_j = s_j s_i$ ($|i-j| \geq 2$) and $s_is_{i+1}s_i = s_{i+1} s_i s_{i+1}$.
Moreover, the abstract group defined by the $s_i$ with the given relations is
a Coxeter group, isomorphic to $S_n$.

This set $\set{s_i}$ can be presented by a graph on the vertices $1,\dots,n$, where $s_i$
is the edge connecting $i$ and $i+1$.
More generally, one can use any connected graph $T$ on $n$ vertices to define a Coxeter group $\Cox{T}$,
from which there is a natural projection onto the corresponding symmetric group. 

The kernel of this projection is generated by elements coming
from two families; one corresponding to triples of vertices in $T$ which
meet in a common vertex, and one to the
cycles of $T$. Let $\CoxY{T}$ denote the quotient of $\Cox{T}$ obtained by assuming the
first family of relations to hold.

One motivation to study the groups $\CoxY{T}$ comes from algebraic
geometry, where these groups are a key ingredient in studying
certain invariants of surfaces (see \cite{T1} for a discussion on
the computation of those invariants). {}From another direction,
signed graphs are used in \cite{CST} to define generalized Coxeter
groups, which are quotients of ordinary Coxeter groups. Our groups
$\CoxY{T}$ belong to this class. The generalized Coxeter groups
which our results enable us to compute are discussed in Subsection
\ref{others}. These include the group $D_2$ whose computation
occupies a large portion of \cite{FJNT}, and a certain family of
Tsaranov groups. For example, the Tsaranov group of a hexagon
(which attracted much attention, see \cite[Example 8.6]{CST}), is
identified in Corollary \ref{TsaranovHexagon}.

There are other indications that $\CoxY{T}$ is a natural quotient
of the Coxeter group $\Cox{T}$. For example, their parabolic
subgroups are well behaved: if $T' \sub T$ is a subgraph, then the
subgroup of $\CoxY{T}$ generated by the elements of $T'$ is
isomorphic to the abstract group defined on $T'$. This is shown in
Subsection \ref{parabolicss}.

We prove that
(with the exception of $T$ equals $K_4$, the complete graph on $4$ vertices),
$\CoxY{T}$ is contained in
the semidirect product $\semidirect{S_n}{\pi_1(T)^n}$, thus
solving the word problem for these groups.
On the other hand $\CoxY{T}$ contains copies of $\pi_1(T)$, showing that it is
virtually solvable (that is, has a solvable subgroup of finite index) iff $T$ has at most one cycle.
Moreover, if $T$ has one cycle then $\CoxY{T}$ is
virtually Abelian. This supports Teicher's conjecture that the invariants mentioned above are either
virtually solvable, or contain a free subgroup \cite{T0}. 
\forget Our identification of $\CoxY{T}$ is used in \cite{TxT} to
show that the fundamental group of the Galois cover of the surface
$\T\times \T$ with respect to a generic projection to $\C\P^2$, is
virtually solvable. \forgotten

Recently, Margulis and Vinberg \cite[Cor. 2]{MV} proved that
infinite non-affine Coxeter groups are large (\ie\ virtually have
a free quotient). In particular the group $\Cox{T}$ is large for
every graph $T$ other than a line, a $Y$-shaped graph (on four
vertices), or a cycle. Our results provide more information on
$\Cox{T}$, proving that already $\CoxY{T}$ is large if $T$ has at
least two cycles, while the kernel of $\Cox{T}\ra \CoxY{T}$ is
large
otherwise. 

In Section \ref{base} we give the basic
 definitions and properties, and
briefly describe an application for our results to algebraic
geometry. Spanning subtrees of $T$ are an important tool
throughout, and in Section \ref{parabolicsec} we prove that the
subgroup generated by a spanning subtree
is the symmetric group. We then describe 
an action
of $S_n$ on the kernel of the projection $\CoxY{T} \ra S_n$, which uses two different embeddings of $S_n$ to $\CoxY{T}$.
In Section \ref{mainsec} we prove the main result,
that this kernel is
isomorphic to a given abstract group, given by generators and relations. This
group is studied in Section \ref{twogps}, where we show it embeds in a direct product of free groups.

The applications to Theorem \ref{main}
are given in Section \ref{appl}: in Corollary \ref{whensolve} we
give the criterion for $\CoxY{T}$ to be virtually solvable. Another immediate
result is that $\CoxY{T}$ depends only on the number of vertices and cycles
of $T$. 
In Subsection \ref{CoxGraph} we discuss the Coxeter graph of $\Cox{T}$
and some special cases.

\section{Presentations of $S_n$ on transpositions}\label{base}

Let $T$ be a graph on $n$ vertices.
Consider the group
generated by the transpositions $(ab) \in S_n$ for the edges $(a,b)$ in
$T$; obviously this is the full symmetric group $S_n$ iff $T$ is
connected. Throughout the paper, all our graphs are simple (\ie\ no repeated
edges or loops).

Recall that a Coxeter group is a group with generators $s_1,\dots,s_k$,
and defining relations $s_i^2 = 1$
and $(s_i s_j)^{m_{ij}} = 1$, where $m_{ij} \in \set{2,3,\dots,\infty}$.
The finite Coxeter groups are completely classified (see \cite{Bourbaki}),
and they are the finite (real) reflection groups.

We use the graph $T$ to define a Coxeter group $\Cox{T}$, as
follows.

\begin{defn}\label{CT}
The group $\Cox{T}$ is generated by the edges $u \in T$, subject to
the following relations:
\begin{equation}\label{square}
u^2=1 \qquad \mbox{for all } u\in T,
\end{equation}
\begin{equation}\label{commute}
uv=vu \qquad \mbox{if } u,v \mbox{ are disjoint, {\emph{and}}}
\end{equation}
\begin{equation}\label{triple}
uvu=vuv \qquad \mbox{if } u,v \mbox{ intersect.}
\end{equation}
\end{defn}
Note that the last relation is equivalent to $(uv)^3 = 1$, so $C =
\Cox{T}$ is indeed a Coxeter group.
\begin{defn}\label{mudef}
Let $T$ be a graph. The map
$$\maptoSn \co \Cox{T} \ra S_n$$
is defined by sending $u = (a,b)$ to the transposition $(ab)$.
\end{defn}
This map is easily seen to be well defined. It is natural to ask
what relations we need to add to $\Cox{T}$ in order to obtain a
presentation of $S_n$, or in other words what elements generate
the kernel of $\maptoSn$. Let $u,v,w \in T$ be three edges meeting
in the vertex $a$, and let $a_u,a_v,a_w$ denote the other vertices
of the respective edges. The transposition $(a,a_u)$ commutes with
$(a,a_v)(a,a_w)(a,a_v) = (a_v,a_w)$. This motivates the following
relation:
\begin{equation}\label{fork}
\, [u,vwv] = 1 \qquad \mbox{for $u,v,w \in T$ which meet in a vertex.}
\end{equation}
Let $u_1,\dots,u_m \in T$ be a cycle. In that we mean that each $u_i$ shares a common vertex with $u_{i-1}$ and $u_{i+1}$
(and $u_1$ with $u_m$), and there are no other intersections. One easily checks that the following relation holds in $S_n$:
\begin{equation}\label{cycle}
\, u_1 \dots u_{m-1} = u_2 \dots u_m.
\end{equation}
It turns out that the relations listed above are enough:
\begin{thm}[{\cite[Prop. 3.4]{Serg}}]\label{Snpres}
Let $T$ be a connected graph on $n$ vertices.
The symmetric group $S_n$ has a presentation with the edges of
$T$ as generators, and the relations \eqs{square}{triple}, \eq{fork} and \eq{cycle}.
\end{thm}

The object we study in this paper is the group $\CoxY{T}$, which we now define.
The subscript ${}_Y$ symbolizes the three edges meeting in a vertex.
\begin{defn}\label{CYT}
Let $T$ be a graph on $n$ vertices.
$\CoxY{T}$ is the group generated by the edges of $T$, with the relations
\eqs{square}{fork}.

The map of Definition \ref{mudef} induces a map $\maptoSn \co
\CoxY{T} \ra S_n$. We let $K(T) \sub \CoxY{T}$ denote the kernel
of this map.
\end{defn}

If $T$ is a tree (\ie\ a connected graph with no cycles), then
\eq{cycle} is vacuously satisfied, so from Theorem \ref{Snpres} we
obtain
\begin{cor}\label{Artree} 
If $T$ is a tree then $\CoxY{T} \isom S_n$.
\end{cor}

The main result of this paper is Theorem \ref{main},
which presents $K(T)$ (defined above)
as a certain subgroup of the direct product $\pi_1(T)^n$. It follows
(Corollary \ref{whensolve}) that $\CoxY{T}$ is virtually solvable iff $T$ has at most one cycle, in which case
it is virtually Abelian.

One application of the study of $\CoxY{T}$ is in algebraic
geometry, specifically to the computation of fundamental groups of
complements of plane 
curves in algebraic surfaces. This problem goes back to Zariski,
and was, in part, the motivation behind van Kampen's celebrated
theorem. In recent years Moishezon and Teicher have made
systematic attempts to compute the fundamental groups of certain
canonical configurations explicitly in terms of generators and
relations, using this to obtain information
on the structure of these groups, see 
\cite{T1}.

In these computations, the first step is to use van Kampen's
theorem to obtain a presentation of the fundamental group. This
group is naturally generated by pairs of elements
$\Gamma_j,\Gamma_{j'}$ ($j = 1,\dots,N$ for some $N$), and usually
has many defining relations. Let $G$ denote the fundamental group,
modulo the relations $\Gamma_j^2 = \Gamma_{j'}^2 = 1$. There is a
natural map $\phi_G$ from $G$ onto a symmetric group $S_n$,
sending the generators $\Gamma_j, \Gamma_{j'}$ to transpositions
(where $\phi_G(\Gamma_j) = \phi_G(\Gamma_{j'})$). It is known that
the fundamental group of a Galois cover of the surface, with
respect to a generic projection to $\C\P^2$, is the kernel of
$\phi_G$. This important invariant of surfaces is useful for
classification of moduli spaces.

One can formally define a quotient group $C= \sg{u_j}$ of $G$ by
applying $\theta \co \Gamma_j,\Gamma_{j'} \mapsto u_j$ to the
relations of $G$. Then the map $\phi_G$ splits as $\phi_G = \phi
\circ \theta$, where $\phi \co C \ra S_n$ is defined by $\phi(u_j)
= \phi_G(\Gamma_j)$. It turns out that $C$ is a Coxeter group of
the form $\Cox{T}$ (for a certain graph $T$ on $n$ vertices), with
some extra relations: relation \eq{cycle} holds for \emph{some} of
the cycles in $T$. As the next remarks show, in the presence of
enough cyclic relations, $C$ actually becomes a quotient of
$\CoxY{T}$, and so our computation of $\CoxY{T}$ allows to compute
$C$.

\forget
 Few remarks concerning the relations \eq{fork} and
\eq{cycle} are in order; then we get back to explain our interest
in the groups $\CoxY{T}$. \forgotten

\begin{rem}\label{forkall}
If $u,v,w \in T$ meet in a vertex, then the relations $[u,vwv] = 1$, $[v,wuw] = 1$ and $[w,uvu] = 1$ are all equivalent (since $[v,wuw] = w[wvw,u]w$).
\end{rem}

\begin{rem}\label{cycleall}
If $u_1,\dots,u_m$ form a cycle in $T$, then all the relations of the form of \eq{cycle}
corresponding to that cycle are equivalent: it does not matter which edge is labelled $u_1$,
nor in which direction the edges are labelled.
\end{rem}
\begin{proof}
Assume $u_1 \dots u_{m-1} = u_2 \dots u_m$. Multiplying by $u_2$ from the left we obtain
$u_1 u_2 u_1 u_3 \dots u_{m-1} = u_2 u_1 u_2 \dots u_{m-1} = u_3 \dots u_m$, but since $u_1$ commutes
with $u_3,\dots,u_{m-1}$, we get $u_1 u_2 \dots u_{m-1} = u_3 \dots u_m u_1$. Thus rotation of
the labels preserves the equality. Taking inverses in \eq{cycle} we see that this is
also the case with changing the direction.
\end{proof}

These two observations justify calling \eq{fork} \emph{the}
relation induced by the triple $u,v,w$, and \eq{cycle} \emph{the}
equation induced by the cycle. We conclude with the following
remark, which explains how cyclic relations imply 'triple'
relations.
\begin{rem}
Let $u_1,\dots,u_m$ be a cycle in $T$,
and $w$ an edge meeting the cycle at the vertex $u_i \cap u_{i+1}$.
Then the relation \eq{fork} corresponding to $u_i,u_{i+1},w$
follows from the relation \eq{cycle} associated to the cycle.
\end{rem}
\begin{proof}
We may assume $i = 1$. Then by assumption we have
$u_2 u_1 u_2 = u_3 \dots u_{m-1} u_m u_{m-1} \dots u_3$, and the $u_i$ ($i \geq 3$) commute with $w$, showing that $[w,u_2 u_1 u_2] = 1$.
\end{proof}

\section{The subgroup of $\CoxY{T}$ generated by a tree}\label{parabolicsec}

In later sections, it will be useful to know that the subgroup of
$\CoxY{T}$ generated by the edges of a spanning subtree $T_0$ of
$T$, is isomorphic to the abstract group defined on $T_0$, \ie\ to
$\CoxY{T_0} = S_n$.

We start with the following easy observation.
\begin{rem}\label{parabolicweak}
Let $T' \sub T$ be any subgraph. Then the subgroup $\sg{T'}$ of $\CoxY{T}$ generated by the
edges of $T'$ is a quotient of $\CoxY{T'}$.
\end{rem}
\begin{proof}
The map $\CoxY{T'} \ra \sg{T'}$ defined by $u \mapsto u$ is well defined,
since every relation in $\CoxY{T'}$ is assumed to hold in $\CoxY{T}$.
Its image is obviously $\sg{T'}$, so we are done.
\end{proof}

Though it is possible to prove directly that the map $\pi \co \CoxY{T'} \ra \sg{T'}$ is an isomorphism, 
we postpone further treatment of the subject to Section \ref{appl}
(Proposition \ref{parabolic}), where it will be derived as an easy consequence
of the main results.

\begin{prop}\label{parabolicSn}
Let $T$ be a graph, and $T_0 \sub T$ a subtree.
Let $\sg{T_0}$ be the
subgroup of $\CoxY{T}$ which is generated by the vertices $u \in T_0$.
Then $\sg{T_0} \isom S_n$.
\end{prop}
\begin{proof}
By the 
remark, $\sg{T_0}$ is the image of the map $\pi \co \CoxY{T_0} \ra
\CoxY{T}$ defined by $u \mapsto u$.

In the diagram below, $\maptoSn, \maptoSn_0$ are the maps of Definition
\ref{mudef} for $T, T_0$, respectively. The diagram commutes by
definition of $\pi$.

\ifXY
\begin{equation}\nonumber
\xymatrix{
\CoxY{T_0} \ar[rd]^{\maptoSn_0} \ar[d]^{\pi}   &  \\
\CoxY{T} \ar[r]^{\maptoSn}  & S_n \\
}
\end{equation}
\else

\bigskip

(A commutative diagram for $\maptoSn,\maptoSn_0$).

\bigskip
\fi 
By Corollary \ref{Artree}, $\maptoSn_{0}$ is an isomorphism, so
the composition $(\maptoSn_0^{-1} \maptoSn) \circ \pi$ is the
identity on $\CoxY{T_0}$. On the other hand if $u\in T_0$ then
$\maptoSn_0^{-1}\maptoSn(u)$ is the generator $u \in \CoxY{T}$, so
that $\pi \circ (\maptoSn_0^{-1} \maptoSn\restricted{\sg{T_0}}) =
1_{\sg{T_0}}$. This proves that $\pi \co \CoxY{T_0} \ra \sg{T_0}$
is an isomorphism.
\end{proof}

\section{$S_n$-action on cycles}\label{actionsec}

In this section we focus on a single cycle in $T$. We study the
cyclic relation \eq{cycle} corresponding to this cycle, and
describe an action of the symmetric group on the normal subgroup
it generates. Recall that $n$ is the number of vertices in $T$.
Throughout, we multiply permutations by $(\s \ta)(a) = \ta
(\s(a))$.

Let $u_1,\dots,u_m$ denote the
edges of the cycle ($m \leq n$ is the length of the cycle),
and renumber the vertices so that
$u_i = (i-1, i)$ ($i = 2,\dots,m$), with $u_1 = (1,m)$.
Let
\begin{equation}\label{gammacycdef}
\gamma_i = u_{i+2}\dots u_m u_1 \dots u_{i}, \qquad i = 1,\dots,
m-2,
\end{equation}
and also set  $\gamma_{m-1} = u_1 \dots u_{m-1}$ and $\gamma_m =
u_2 \dots u_m$.

According to Theorem \ref{Snpres}, the kernel of $\maptoSn
:\CoxY{T} \ra S_n$ is the normal subgroup generated by elements of
the form $\gamma_{m-1}^{-1}\gamma_m$ (ranging over all the cycles
of $T$). In Remark \ref{cycleall} we have seen that $\gamma_j^{-1}
\gamma_i$ belongs to this kernel for all $i,j = 1, \dots, m$.
Choose a spanning subtree $T_0$ of $T$. It would be nice to have a
natural action of the subgroup $S_n  = \sg{T_0}$ on the set
$\set{\gamma_i^{-1}\gamma_j}$; such an action can take the form
\begin{equation}\label{doubleact}
\s^{-1} \gamma_j^{-1} \gamma_i \s = \gamma_{\s j}^{-1}\gamma_{\s i}, \qquad \s \in \sg{T_0}.
\end{equation}

For \eq{doubleact} to even make sense,
we need $\gamma_i$ to be defined for every
$i = 1,\dots,n$, not just for vertices on the cycle.
We will show that assuming \eq{doubleact} to hold, we are led to a unique
definition of the $\gamma_i$. We will later show that using this
definition, Equation \eq{doubleact} is indeed
satisfied.
Since $\gamma_j^{-1} \gamma_i$ is mapped to the unit element of $S_n$, all the $\gamma_i$ map to the same
permutation $\maptoSn(\gamma_1) = (23)\dots(m-1\,m)(1m) = (m \dots 3 2 1)$.

A \defin{path} in a graph is an ordered list of edges, in which every two
consecutive edges have one vertex in common, and there are no other
intersections.

\begin{prop}\label{theydo}
Suppose that elements $\gamma_i \in \CoxY{T}$ ($i = 1,\dots,n$) are defined for every cycle in $T$,
such that \Eq{doubleact} holds, and for the vertices $i$ on the cycle,
$\gamma_i$ is defined by \Eq{gammacycdef}.

Then $K(T) = \Ker(\maptoSn \co \CoxY{T}\ra S_n)$ is generated by
the elements $\gamma_j^{-1}\gamma_i$ as a subgroup (rather than a
normal subgroup) of $\CoxY{T}$.
\end{prop}
\begin{proof}
Let $N$ denote the subgroup of $\CoxY{T}$ generated by the
elements $\gamma_j^{-1}\gamma_i$ ($i,j = 1,\dots,n$) for all
cycles in $T$. We already know that $K(T)$ is generated by the
$\gamma_{m-1}^{-1}\gamma_m$ as a normal subgroup. \Eq{doubleact}
indicates that the $\gamma_j^{-1}\gamma_i$ are all conjugate to
each other, so $N \sub K(T)$.

We need to show that $N$ is normal. Let $\gamma_b^{-1}\gamma_a$ be
one of the generators, and let $x \in T$. If $x \in T_0$, then by
\Eq{doubleact}, $x \gamma_b^{-1} \gamma_a x^{-1}$ is of the same
form. Otherwise, there is a unique path in $T_0$ connecting the
two vertices of $x$, and together with $x$ this is a cycle in $T$.
Denote by $\delta_{i'}$ the elements defined in \eq{gammacycdef}
for that cycle (on the vertices $1',\dots,m'$), where $x$ is the
edge connecting $1'$ and $2'$. Then  $x \delta_{1'}^{-1}
\delta_{m'} \in \sg{T_0}$. Let $\s = \maptoSn(x\delta_{1'}^{-1}
\delta_{m'})$. Then $x \gamma_b^{-1} \gamma_a x^{-1} = \s
(\delta_{m'}^{-1}\delta_{1'})(\gamma_b^{-1}
\gamma_a)(\delta_{1'}^{-1}\delta_{m'}) \s^{-1}$ which is in the
subgroup generated by the $\gamma_j^{-1}\gamma_i$ and the
$\delta_j^{-1}\delta_i$. 
\end{proof}

{}From Equation \eq{doubleact} alone it is not clear how a
conjugate of $\gamma_i$ by $\s$ looks like. We cannot expect
$\s^{-1} \gamma_i \s$ to be of the form $\gamma_j$, since
$\maptoSn(\gamma_i)$ is not a central element of $S_n$. Still,
from Equation \eq{doubleact} it follows that $\gamma_i \s
\gamma_{\s i}^{-1}$ is independent of $i$, so for every $\s \in
\sg{T_0}$ there is some $\til{\s} \in \CoxY{T}$ such that
\begin{equation}\label{singleact}
\til{\s}^{-1} \gamma_a \s = \gamma_{\s a}, \qquad a = 1,\dots,n.
\end{equation}

{}From this it would follow that $\tilmap[\s]$ is a homomorphism.
We will show that this is the case after the $\til{\s}$ are
specified.

Equation \eq{singleact} can serve as the definition of $\til{\s} = \gamma_a \s \gamma_{\s a}^{-1}$
(ignoring for the moment the fact that $\gamma_a$ is not always defined, and that we need to show
the definition is independent of $a$). Fixing the cycle, let us compute $\gamma_a u \gamma_{u(a)}^{-1}$
for arbitrary $u \in T_0$ (and $a$ of our choice).
There are four cases to consider.
\begin{enumerate}
\item $u$ does not touch the cycle. Then choosing $a$ on the cycle we have that $\gamma_a$ commutes with $u$, so that
$\gamma_a u \gamma_a^{-1} = u$.
\item $u$ is part of the cycle, say $u = u_i$ ($i = 1,\dots,m$). Taking $a = i-2$ (with $a = m$ if $i = 2$ and $a = m-1$ if $i = 1$),
we see that
\begin{eqnarray*}
\gamma_{i-2} u_i \gamma_{u(i-2)}^{-1}
    & = & \gamma_{i-2} u_i \gamma_{i-2}^{-1} \\
    & = & u_{i}\dots u_m u_1 \dots u_{i-2} u_i u_{i-2}\dots u_1 u_m \dots u_{i} \\
    & = & u_{i} u_{i+1} u_i u_{i+1} u_{i} \\
    & = & u_{i+1},
\end{eqnarray*}
where for $i = m$ the calculation gives $u_1$.
\item $u$ touches the cycle at one vertex, $i$ ($1 \leq i \leq
m$). Taking $a = i-1$, we have
\begin{eqnarray*}
\gamma_{i-1} u \gamma_{u(i-1)}^{-1}
    & = & \gamma_{i-1} u \gamma_{i-1}^{-1} \\
    & = & u_{i+1}\dots u_m u_1 \dots u_{i-1} u u_{i-1}\dots u_1 u_m \dots u_{i+1} \\
    & = & u_{i+1} u u_{i+1}.
\end{eqnarray*}
Notice that if $\phi(u) = (i,k)$, then $\phi(u_{i+1} u u_{i+1}) =
(i+1,k)$.
\item \label{touchtwice}
$u$ touches the cycle at two vertices, $i,j$, with $j>i+1$. Then $u(i-1) = i-1$ and
\begin{eqnarray*}
\gamma_{i-1} u \gamma_{i-1}^{-1}
    & = & u_{i+1}\dots u_m u_1 \dots u_{i-1} u u_{i-1}\dots u_1 u_m \dots u_{i+1} \\
    & = & u_{i+1}\dots u_{j+1} u u_{j+1} \dots u_{i+1} \\
    & \eqY & u_{i+1} u_{j+1} u u_{j+1} u_{i+1}.
\end{eqnarray*}
Here $\phi(u) = (i,j)$ by assumption, and $\phi(u_{i+1} u_{j+1} u
u_{j+1} u_{i+1}) = (i+1,j+1)$.
\end{enumerate}

We take the results of this computation as our definition.
\begin{defn}\label{tildedef}
Fixing a (directed) cycle in $T$, for every $u \in T$ we define $\til{u} \in
\CoxY{T}$, as follows:
$$\til{u} = \begin{cases}
u & \mbox{if $u$ does not touch the cycle} \\
u_{i+1} & \mbox{if $u = u_i$ on the cycle} \\
u_{i+1} u u_{i+1} & \mbox{if $u$ touches the cycle at vertex $i$ only} \\
u_{i+1} u_{j+1} u u_{j+1} u_{i+1} & \mbox{if $u$ touches the cycle at the vertices $i,j$}. \\
\end{cases} $$
\end{defn}

\begin{rem}\label{actil}
Let $\tau = (m \dots 3 2 1)$. {}From the definition it follows
that for every $u \in T$, $\phi(\til{u}) = \tau \phi(u)
\tau^{-1}$.
\end{rem}

{}From now and to the rest of the paper, whenever we enumerate the
edges of a cycle $u_1,\dots,u_m$, we will choose the enumeration
so that $u_2,\dots,u_m \in T_0$ (so that necessarily $u_1 \not \in
T_0$). Under this assumption, an edge $u \in T_0$ cannot touch the
cycle more than once. As we shall see in Section \ref{mainsec},
this makes our cycles "basic cycles" (with respect to $T_0$).

\begin{cor}\label{ishom}
The map $\tilmap[u]$ extends to a homomorphism $\tilmap[\s]$ from
$\sg{T_0}$ to $\CoxY{T}$.
\end{cor}
\begin{proof}
We only need to show that $\tilmap[u]$ preserves the defining
relations of $\sg{T_0} \isom \CoxY{T_0}$. Let $T_0'$ be an
arbitrary spanning subtree of $T$. If $R(u)$ is a relator on the
generators of $\CoxY{T_0}$, then $\phi(R(\til{u})) =
R(\phi(\til{u})) = R(\tau \phi(u) \tau^{-1}) = \tau \phi(R(u)) \tau^{-1} = 1$,
showing that $R(\til{u}) \in \Ker(\phi) = K(T)$. Thus, if for all
the generators $u$ involved in a relator we have $\til{u} \in
\sg{T_0'}$, then $R(\til{u}) \in \sg{T_0'} \cap K(T) = 1$.

Let $u,v \in T_0$; then the subgraph of $T$ built of the
generators participating in $\til{u},\til{v}$ never contains a
cycle. So it can be completed to a spanning subtree $T_0'$,
showing that the relations in which $u,v$ are involved are
satisfied by $\til{u},\til{v}$. The same argument applies to
$u,v,w \in T_0$ which meet in a common vertex, showing that
relation \eq{fork} is also preserved, so we are done.
\end{proof}

In Subsection \ref{tildess} we will improve this result, and show
that $\tilmap[u]$ actually extends to an automorphism of
$\CoxY{T}$. However it should be emphasized that in general
\Eq{singleact} does not hold for every $u \in T$.

Now that $\til{\s}$ is defined, we can use \eq{singleact} to
define $\gamma_a$ for arbitrary $a$: let $v_1,\dots,v_s$ denote
the path in $T_0$ from the vertex $1$ to $a$, and define
\begin{equation}\label{gammadef}
\gamma_{a} = \til{v_s} \dots \til{v_1} \gamma_1 v_1 \dots v_s. 
\end{equation}

We first need to check that this new
definition coincides with \eq{gammacycdef}.
Indeed, if $a$ is on the cycle then the path from $1$ to $a$
is $u_2,\dots,u_{a}$, so
\begin{eqnarray*}
\til{u_a} \dots \til{u_2} \gamma_1 u_2 \dots u_a & = &
        u_{a+1} \dots u_3  \cdot u_{3}\dots u_m u_1 \cdot  u_2 \dots u_a \\
    & = & u_{a+2} \dots u_m u_1 u_2 \dots u_a
\end{eqnarray*}
which is the previous definition.
Secondly, if $b$ is a vertex on the path from $1$ to $a$,
say $b = v_{k-1} \cap v_{k}$, then by definition of $\gamma_b$,
$\gamma_a = \til{v_s} \dots \til{v_k} \gamma_b v_k \dots v_s$.

\begin{rem}\label{taueverywhere}
Recall that $\tau = \phi(\gamma_1) = (m \dots 3 2 1)$. For every
$a$ we have by Remark \ref{actil} that $\phi(\gamma_a) = \tau
\phi(v_s \dots v_1) \tau^{-1} \phi(\gamma_1) \phi(v_1 \dots v_s) =
\tau$.
\end{rem}

\begin{prop}\label{singleactprop}
With the above definitions 
of the $\gamma_i$ and $\til{u}$, Equation \eq{singleact} holds for
every $\s \in \sg{T_0}$ and $a = 1, \dots, n$.
\end{prop}
\begin{proof}
It is enough to assume $\s = u \in T_0$, and we need to check that
$\til{u}^{-1} \gamma_a u = \gamma_{u(a)}$. We use induction on the distance $s$ of $a$ from the vertex~$1$.
First assume that $a = 1$. There are four cases to consider.
        \begin{enumerate}
            \item $u$ touches the vertex $1$. Then by definition $\gamma_{u(1)} = \til{u} \gamma_1 u$.
            \item $u$ is on the cycle, $u = u_i$ with $2 < i \leq m$. Then
\begin{eqnarray*}
\gamma_1 u \gamma_{u(1)}^{-1}
& = & u_{3} \dots u_{m} u_1 u_{i} u_1 u_m \dots u_{3} \\
& = & u_{3} \dots u_i u_{i+1} u_{i} u_{i+1} u_i \dots u_{3} \\
& = & u_{3} \dots u_{i-1} u_{i+1} u_{i-1} \dots u_{3} \\
& = & u_{i+1} = \til{u}.
\end{eqnarray*}
            \item $u$ touches the cycle at a vertex $i$ other than $1$. Then $\til{u} = u_{i+1} u u_{i+1}$, and
\begin{eqnarray*}
\gamma_1 u \gamma_{u(1)}^{-1}
& = & u_{3} \dots u_m u_1 u u_1 u_m \dots u_{3} \\
& = & u_{3} \dots u_{i+1} u u_{i+1} \dots u_{3} \\
& \eqY & u_{i+1} u u_{i+1} = \til{u}.
\end{eqnarray*}
            \item $u$ does not touch the cycle. Then $\til{u} = u$ commutes with
            every $u_i$ on the cycle, so $\til{u}^{-1} \gamma_1 u  = \gamma_1$.
        \end{enumerate}

Now suppose $a \neq 1$, and let $v_1,\dots,v_s$ denote the path
from $1$ to $a$. Let $b = v_s(a)$; then $v_1,\dots,v_{s-1}$ is the path from $1$ to $b$.
By the definition of $\gamma_b$ and $\gamma_a$, and the induction hypothesis, we have that
\begin{eqnarray*}
\til{u} \gamma_a u
    & = & \til{u} \til{v_s} \gamma_b v_s u \\
    & = & \til{u} \til{v_s} \til{u} \cdot \til{u} \gamma_b u \cdot u v_s u \\
    & = & \til{u} \til{v_s} \til{u} \cdot \gamma_{u(b)} \cdot u v_s u.
\end{eqnarray*}
In particular, if $u$ is disjoint from $v_s$ then
$\til{u}\til{v_s}\til{u} = \til{v_s}$ by \Cref{ishom}, $u(b) = b$,
and $\til{u} \gamma_a u  = \til{v_s} \gamma_{b} v_s = \gamma_a$.
Thus, we may assume that $u$ touches $a$ or $b$. There are two
cases to check if $u$ touches $a$: either $u = v_s$, in which case
$u(a) = b$ and $\til{u}^{-1} \gamma_a u = \til{v_s} \gamma_a v_s =
\gamma_b = \gamma_{u(a)}$; or $u(a) \neq b$, where $\gamma_{u(a)}
= \til{u} \gamma_a u$ by definition of $\gamma_{u(a)}$. Next,
assume $u(a) = a$, but $c = u(b) \neq b$.

If $u = v_{s-1}$, let $b' = v_{s-1}(b)$. Then by the induction hypothesis
$\til{u} \gamma_a u = \til{u} \til{v_s} \gamma_b v_s u = \til{v_s} \til{u} \til{v_s} \til{u} \gamma_b u v_s u v_s = \til{v_s} \til{u} \til{v_s} \gamma_{b'} v_s u v_s = \til{v_s} \til{u} \gamma_{b'} u v_s = \til{v_s} \gamma_{b} v_s = \gamma_a$.
So we may assume $u$ is not part of the path $v_1,\dots,v_s$; see Figure \ref{pic1}.
We need to show that
$\til{u} \gamma_a u = \gamma_{a}$; writing $\gamma_a = \til{v_s} \gamma_b v_s$,
this is equivalent to $\til{v_s} \til{u} \til{v_s} \gamma_b v_s u v_s = \gamma_b$, which
is symmetric under replacing $a$ and $v_s$ by $c$ and $u$.

\begin{figure}[!h]
\ifXY
\begin{equation}\nonumber
\xymatrix{ 
{}
    &
    &
    &
    &
    & c
    &   \\
{}
    & 1 \ar@{-}[r]^{v_1} \ar@{=}@/_4pt/[ld]^{u_1}
    &   \ar@{.}[rr]
    &
    &   \ar@{-}[r]^{v_{s-1}}
    & b \ar@{-}[r]^{v_s} \ar@{-}[u]^{u}
    & a \\
m \ar@{-}[d]^{u_m}
    &
    &
    &
    &
    &
    &   \\
m-1 \ar@{.}@(d,r)@/_14pt/ [rrd] 
    &
    &
    &
    &
    &
    &   \\
{}
    &
    & \ar@{.}@(r,l)@/_42pt/ [uuu]
    &
    &
    &
    &   \\
}
\end{equation}
\else

\bigskip

(A picture of the circle).

\bigskip
\fi 
\caption{}\label{pic1}
\end{figure}

There are three cases to consider.
\begin{enumerate}
    \item $a$ is on the cycle. Then $b$ is also on the cycle, but $c$ is not.
Compute:
\begin{eqnarray*}
\gamma_a u \gamma_{u(a)}^{-1}
& = & u_{a+2} \dots u_m u_1 \dots u_a u u_a \dots u_1 u_m \dots u_{a+2} \\
& \eqY & u_a u u_a  = \til{u}.
\end{eqnarray*}
    \item $a,c$ are not on the cycle but $b$ is. Then
\begin{eqnarray*}
\gamma_a u \gamma_{u(a)}^{-1}
    & = & \til{v_s} \gamma_b v_s u v_s \gamma_{b}^{-1} \til{v_s} \\
    & \eqY & \til{v_s} v_s u v_s \til{v_s} \\
    & = & u_{b+1} v_s u_{b+1} v_s u v_s \cdot u_{b+1} v_s u_{b+1} \\
    & \eqY & u_{b+1} v_s v_s u v_s v_s u_{b+1} = u_{b+1} u u_{b+1} = \til{u}.
\end{eqnarray*}
    \item $b$ is not on the cycle. Let $w_1,\dots,w_p$ be the path from a
vertex $j$ on the cycle to $b$ (so that $w_p = v_{s-1}$).
\begin{eqnarray*}
\gamma_a u \gamma_{u(a)}^{-1}
    & = & \til{v_s} \til{w_p} \dots \til{w_1} \gamma_j w_1 \dots w_p v_s u v_s w_p \dots w_1 \gamma_j^{-1} \til{w_1} \dots \til{w_p} \til{v_s} \\
    & \eqY & \til{v_s} \til{w_p} \dots \til{w_1} \gamma_j v_s u v_s \gamma_j^{-1} \til{w_1} \dots \til{w_p}\til{v_s} \\
    & \eqY & \til{v_s} \til{w_p} v_s u v_s \til{w_p}\til{v_s} \\
    & = & v_s w_p v_s u v_s w_p v_s \\
    & \eqY & v_s v_s u v_s v_s = u = \til{u}.
\end{eqnarray*}
\end{enumerate}

\end{proof}

Multiplying $\til{\s}^{-1} \gamma_i \s = \gamma_{\s i}$ and the
inverse $\s^{-1} \gamma_j^{-1} \til{\s} = \gamma_{\s j}^{-1}$, we
now get
\begin{cor}\label{PsiPreserveAct}
Equation \eq{doubleact} holds for every $\s \in \sg{T_0}$.
\end{cor}

We conclude with the following properties of the $\gamma_i$.
\begin{prop}\label{gammagammam}
Let $\tau = \phi(\gamma_1)$. For every $j = 1,\dots, n$, we have
that
\begin{equation}\label{ggm}
\gamma_{j} \gamma_{m} = \gamma_{1} \gamma_{\tau(j)}.
\end{equation}
\end{prop}
\begin{proof}
Since $\gamma_m = u_2 \dots u_m \in \sg{T_0}$, we have that
$\til{\gamma_m} = \til{u_2} \dots \til{u_m} = u_3 \dots u_m u_1 =
\gamma_1$. Set $\s = \gamma_m$ in Equation \eq{singleact}, to get
$\gamma_1^{-1} \gamma_j \gamma_m = \gamma_{\tau(j)}$ for arbitrary
$j$.
\end{proof}

Similarly, we have
\begin{prop}\label{gammagamma}
For every $i,j$ on the cycle, we have
\begin{equation}\label{gg}
\gamma_{j} \gamma_{i-1} = \gamma_{i} \gamma_{j-1}.
\end{equation}
\end{prop}
\begin{proof}
The equality holds for $i = 1$ by the last proposition. But this
is a relation on the generators $u_1,\dots,u_m$, so we are done by
renaming the edges (it can also be proved by direct computation).
\end{proof}

{}From Proposition \ref{Phigamma} it follows that $\gamma_{j}
\gamma_{\tau(i)} = \gamma_{i} \gamma_{\tau(j)}$  for any $i,j$,
but we will not need this here.

\section{Two abstract groups}\label{twogps}

In this section we describe two abstract groups, $A_{t,n}$ and
$F_{t,n}$, and show that for all the interesting pairs $(t,n)$
(with one exception), they are isomorphic. In the next section we
will prove that the kernel $\Ker(\maptoSn \co \CoxY{T} \ra S_n)$,
denoted by $K(T)$, is isomorphic to $A_{t,n}$ (for a certain $t$
depending on $T$). The isomorphism from $K(T)$ to $A_{t,n}$ is
rather natural, but $A_{t,n}$ is better understood through its
identification with $F_{t,n}$ (which is an explicit subgroup of
the direct power of a free group). For reasons that will later
become clear, if $n = 4$ we are only interested in the case $t
\leq 3$, and if $n < 4$ in the case $t \leq 1$. The definitions
are, however, general.

\subsection{The group $A_{t,n}$}\label{Atnss}

Let $t \geq 0$, $n\geq 1$.
Fix a set $X = \set{x,y,z,\dots}$ of size $t$.
\begin{defn}\label{AtnDef}
The group $A_{t,n}$ is generated by the $n^2\power{X}$ elements
$x_{ij}$ ($x \in X$, $i, j = 1, \dots,n$), with the defining
relations
\begin{eqnarray}
x_{ii} & = & 1, \label{unit} \\
x_{ij}x_{jk} & = & x_{ik}, \label{transitive}  \\ 
x_{jk}x_{ij} & = & x_{ik}, \label{transitivecomm}
\end{eqnarray}
and
\begin{equation}\label{disjoint}
\, [x_{ij},y_{kl}] = 1 \qquad \mbox{if\ $i,j,k,l$ are different}.
\end{equation}
\end{defn}

We will sometimes
use $A_{X,n}$ to specify the set of generators of $A$; this is
the same group as $A_{t,n}$ for $t = |X|$.

\begin{rem}\label{xiscomm}
For every fixed $x \in X$, the subgroup $\sg{x_{ij}} < A_{t,n}$ is
commutative.
\end{rem}
\begin{proof}
First, $x_{ji} = x_{ij}^{-1}$ by \eq{unit} and \eq{transitive}.
For every $i,j$ we have $x_{1i}x_{j1} = x_{ji} = x_{j1}x_{1i}$ by
\eqs[,]{transitive}{transitivecomm}, so $x_{1j}x_{1i} = x_{1i}x_{1j}$ and $\sg{x_{1i}}$ is
commutative. Now $x_{ij} = x_{1i}^{-1}x_{1j}$, and we are done.
\end{proof}

\begin{exmpl}\label{A1n}
For $t = 1$, $A_{1,n}$ is generated by a single set $\set{x_{ij}}$ and is commutative by Remark \ref{xiscomm}.
Relation \eq{disjoint} is then redundant, so $x_{i,i+1}$ are free generators, and $A_{1,n} \isom \Z^{n-1}$.
\end{exmpl}

\begin{prop}\label{xyzredun}
If $n \geq 5$, or if $t \leq 2$, then $A_{t,n}$ satisfies the following relation:
\begin{equation}\label{xyzeq}
 [x_{ij}, z_{ik} y_{kl} z_{ki}] = 1 \qquad \mbox{if\ $i,j,k,l$ are different}.
\end{equation}
\end{prop}
\begin{proof}
First assume $n \geq 5$, so we can choose $r \neq i,j,k,l$. Using the decomposition $z_{ik} = z_{ir} z_{rk}$ and
relation \eq{disjoint}, we compute
\begin{eqnarray*}
\, [x_{ij}, z_{ik} y_{kl} z_{ki}]
    & = & x_{ij} z_{ik} y_{kl} z_{ki} x_{ji} z_{ik} y_{lk} z_{ki}
    \\
    & = & x_{ij} z_{ir} z_{rk} y_{kl} z_{ki} x_{ji} z_{ir} z_{rk} y_{lk} z_{ki}
    \\
    & = & z_{rk} x_{ij} y_{kl} z_{ir} z_{ki} z_{rk} x_{ji} y_{lk} z_{ir} z_{ki}
    \\
    & = & z_{rk} x_{ij} y_{kl} x_{ji} y_{lk} z_{kr} = 1
\end{eqnarray*}

Now suppose $t = |X|  \leq 2$. If in \eq{xyzeq} we set
$z = y$ or $z = x$ the assertion is trivial, so assume $x = y$,
and compute: 
\begin{eqnarray*}
\, [x_{ij}, z_{ik} x_{kl} z_{ki}]
    & = & x_{ij} z_{ik} x_{kl} z_{ki} x_{ji} z_{ik} x_{lk} z_{ki}
    \\
    & = & x_{ij} z_{ik} x_{kj} x_{jl} z_{kl} z_{li} x_{ji} z_{ik} x_{lk} z_{ki}
    \\
    & = & x_{ij} x_{jl} z_{ik} z_{li} x_{kj} x_{ji} z_{kl} z_{ik} x_{lk} z_{ki}
    \\
    & = & x_{il} z_{lk} x_{kj} x_{ji} z_{il} x_{lk} z_{ki}
    \\
    & = & x_{il} x_{ji} z_{lk} z_{il} x_{kj} x_{lk} z_{ki}
    \\
    & = & x_{jl} z_{ik} x_{lj} z_{ki} = 1
\end{eqnarray*}
\end{proof}

Essentially the same proof shows that under the same assumptions,
\begin{equation}\label{xyzeqinn}
 [x_{ij}, z_{ki} y_{kl} z_{ik}] = 1 \qquad \mbox{if\ $i,j,k,l$ are different}.
\end{equation}
also holds in $A_{t,n}$.

The next proposition will be useful in Section \ref{appl}.
\begin{prop}\label{At1At}
$A_{t-1,n}$ is a subgroup (in fact a retract) of $A_{t,n}$.
\end{prop}
\begin{proof}
Let $X' \subset X$ be sets of sizes $t-1$ and $t$, respectively,
and let $w$ denote the element in $X-X'$. Define $\rho' \co
A_{X',n} \ra A_{X,n}$ by $x_{ij} \mapsto x_{ij}$ for $x \in X'$,
and $\rho \co A_{X,n} \ra A_{X',n}$ by $\rho(x_{ij}) = x_{ij}$ for
$x \in X'$, and $\rho(w_{ij}) = 1$. Since all the relations are
preserved, both maps are well defined.

For every generator $x_{ij}$ of $A_{X',n}$ ($x \in X'$, $1\leq i,j
\leq n$), the maps satisfy $\rho \rho'(x_{ij}) = \rho(x_{ij}) =
x_{ij}$, so $\rho \rho'$ is the identity on $1_{A_{X',n}}$. In
particular $\rho'$ is one-to-one, injecting $A_{X',n}$ onto its
image, the subgroup  $\sg{x_{ij} \suchthat x\in X'}$ of $A_{X,n}$.
\end{proof}

\subsection{The group $F_{t,n}$}\label{Ftnss}

Again let $X$ be a set of size $t \geq 0$.
For every $i = 1,\dots, n$, let $F^{(i)}$ denote the free group on
the $t$ letters $x_i$ ($x \in X$). Set $\Fstar = F^{(1)}\times
\dots \times F^{(n)}$. Then $\Fstar$ has generators $x_{i}$ ($x
\in X$, $i = 1,\dots,n$) and defining relations
\begin{eqnarray}
\, [x_{i}, y_j] & = & 1 \qquad (x,y \in X,\  i\neq j) \label{xwithy}.
\end{eqnarray}

The Abelianization map $\ab: F^{(i)} \ra \Z^t$ defined on every
component by $\ab(x_i) = x$ (where $\Z^t$ is thought of as the
free Abelian group generated by $X$), can be extended to a map
$\Fstar \ra\Z^t$ by summing the entries:
\begin{equation}\label{abdef}
\ab(p_1,\dots,p_n) = \ab(p_1)+\dots+\ab(p_n).
\end{equation}
\begin{defn}\label{FtnDef}
$F_{t,n}$ is the kernel of $\ab \co \Fstar \ra \Z^t$.
\end{defn}

\begin{thm}\label{AisomF}
$A_{t,n} \isom F_{t,n}$ in the following cases: $t \leq 1$, or $n \geq 5$,
or $n = 4$ and $t = 2$.

Moreover, if $n = 4$, then $F_{t,n}$ is a quotient of $A_{t,n}$.
\end{thm}
\begin{proof}
If $t = 0$ then both groups are trivial. 
If $t = 1$ then $A_{t,n} = \Z^{n-1}$,
and $\Fstar = \Z^n$ by definition, so $F_{t,n}=\Z^{n-1}$.
We may thus assume $n \geq 4$.

Define a group $\Astar$ with generators $x_{ij}$ and $s_x$ ($x \in
X$, $i,j = 1,\dots, n$), and defining relations
\eqs{unit}{disjoint} and the following:
\begin{eqnarray}
\, [s_x, y_{ij}] & = & [x_{nk},y_{ij}] \qquad (x, y \in X,\ k \neq i,j), \label{sony} \\
\, [s_x, s_y] & = & [x_{ni},y_{nj}] \qquad (i\neq j,\ i,j \neq n).
\label{sons}
\end{eqnarray}

Define a map $\mu \co \Astar \ra \Fstar$ as follows:
\begin{eqnarray*}
\mu(x_{ij}) &=& x_{j}^{-1}x_{i}, \\
\mu(s_x) &=& x_{n}.
\end{eqnarray*}

A routine check shows that $\mu$ respects all the relations
\eq{unit}, \eq{transitive}, \eq{transitivecomm}, \eq{disjoint}, 
\eq{sony} and \eq{sons}, so is well defined.

Also, $\mu(s_x x_{in}) = x_i$, so $\mu$ is onto $\Fstar$. Since
all the relations of $A_{t,n}$ are assumed to hold in $\Astar$,
$x_{ij} \mapsto x_{ij}$ defines a map $\rho \co A_{t,n} \ra
\Astar$. Also define $\ab \co \Astar \ra \Z^{t}$ by $\ab(x_{ij}) =
1$ and $\ab(s_x) = x$.

Now define $\mu' \co \Fstar \ra \Astar$ by
$$ \mu'(x_i) = s_x x_{in}$$
(in particular we have $\mu'(x_n) = s_x x_{nn} = s_x$). The fact
that \eq{xwithy} is preserved by $\mu'$ for $x=y$, follows from
Remark \ref{xiscomm}. If $x\neq y$, choose $k \neq j,n$ and $l
\neq i,k,n$ (this requires $n\geq 4$) and compute:
\begin{eqnarray*}
[s_x x_{in}, s_y y_{jn}]
    & = & s_x x_{in} s_y y_{jn} x_{ni} s_x^{-1} y_{nj} s_y^{-1} \\
    & = & x_{in} s_x y_{jn} s_y s_x^{-1} x_{ni} s_y^{-1} y_{nj}  \\
    & \stackrel{\eq{sony}}{=} & x_{in} x_{nk} y_{jn} x_{kn} s_x s_y s_x^{-1} x_{ni} s_y^{-1} y_{nj}  \\
    & \stackrel{\eq{sony}}{=} & x_{in} x_{nk} y_{jn} x_{kn} s_x s_y s_x^{-1} s_y^{-1} y_{nl} x_{ni} y_{ln} y_{nj}  \\
    & \stackrel{\eq{sons}}{=} & x_{in} x_{nk} y_{jn} x_{kn} x_{nk} y_{nl} x_{kn} y_{ln} y_{nl} x_{ni} y_{ln} y_{nj}  \\
    & \stackrel{\eq{unit}-\eq{transitivecomm}}{=} &  x_{ki} y_{lj} x_{ik} y_{jl}  \\
    & \stackrel{\eq{disjoint}}{=} & 1.
\end{eqnarray*}
This proves that $\mu'$ is well defined. Moreover, it is easy to check
that $\mu \mu' = 1_{\Fstar}$ and $\mu' \mu = 1_{\Astar}$.
We also use $\mu$ to denote the map $A_{t,n} \ra F_{t,n}$
defined by $\mu(x_{ij}) = x_j^{-1}x_{i}$. Then the following diagram commutes:

\ifXY
\begin{equation}\label{AFstar}\nonumber
\xymatrix{
{}
    & A_{t,n} \ar[r]^{\rho} \ar[d]^{\mu}
    & \Astar \ar[r]^{\ \ab} \ar@<-1ex>[d]_{\mu}
        & \Z^t \ar[r]\ar@{=}[d]
        & 1 \\
1 \ar[r]
    & F_{t,n} \ar@{^{(}->}[r]
        & \Fstar \ar[r]^{\ \ab} \ar@<-1ex>[u]_{\mu'}
        & \Z^t \ar[r]
        & 1 \\
}
\end{equation}
\else

\bigskip

(A commutative diagram for $\mu,\mu'$).

\bigskip
\fi 

It is easy to see that as a normal subgroup of $\Fstar$, $F_{t,n}$
is generated by the elements $x_i^{-1} x_n$ and $[x_n,y_n]$
($x,y\in X$). Thus, applying $\mu'$, we see that $\mu'(F_{t,n})$ is
generated by the elements $x_{ni}$ and $[s_x,s_y] =
[x_{ni},y_{nj}]$; but these are the generators of $\Im(\rho)$, showing that $F_{t,n}$ and $\Im(\rho)$ are
isomorphic.

 From this it follows that $F_{t,n}$ is a quotient group of $A_{t,n}$.
In the next proposition we show that $\rho \co A_{t,n} \ra \Astar$
is injective if $n \geq 5$ or $t \leq 2$, thus completing the
proof.
\end{proof}

We now show that if $n \geq 5$ or $t \leq 2$, then the mapping
$\rho \co A_{t,n} \ra \Astar$ defined by $\rho(x_{ij}) = x_{ij}$
is an embedding. Fix an arbitrary order on $X$.  For $z \in X$,
let $\Astar[z]$ denote the group generated by $A_{t,n}$ and the
elements $s_x$ ($x \leq z$), subject to the relations
\eqs{sony}{sons}. Of course, when defining $\Astar[z]$, the range
of the variables $x$ and $y$ in these relations is $x,y < z$
(rather than all $X$). Thus, if $z$ is the maximal element of $X$,
we have that $\Astar[z] = \Astar$.

\begin{prop}\label{AintoAstar}
Assume $n\geq 5$  or $t \leq 2$.

If $u<z$ are consecutive elements of $X$, then $\Astar[z]$ is a
semidirect product of $\Astar[u]$ and $\Z = \sg{s_z}$. In
particular, $A_{t,n}$ is the subgroup of $\Astar$ generated by the
$x_{ij}$.
\end{prop}
\begin{proof}
Let us define an automorphism of $\Astar[u]$ by
\begin{eqnarray*}
\bar{s_z}(x_{ij}) & = & z_{nk}x_{ij}z_{nk}^{-1} \qquad (x \in X, \ k \neq i,j), \\
\bar{s_z}(s_x) & = & [z_{ni},x_{nj}]s_x \qquad (x< z,\ i\neq j,\
i,j \neq n).
\end{eqnarray*}
We must first verify that the definition does not depend on the
choice of $k$ (in the first case), or $i,j$ (in the second). The
equality $z_{nk}x_{ij}z_{nk}^{-1} = z_{nk'}x_{ij}z_{nk'}^{-1}$
(for $k,k' \neq i,j$) follows from \eq{disjoint}. Now suppose that
$i,j,n$ are different, and that $i',j',n$ are different. If $i'
\neq j$ then it is easy to check that $[z_{ni},x_{nj}] =
[z_{ni'},x_{nj}] = [z_{ni'},x_{nj'}]$; likewise the equality
follows if $i \neq j'$. In order to show that $[z_{ni},x_{nj}] =
[z_{nj},x_{ni}]$, choose $k \neq i,j,n$; then $[z_{ni},x_{nj}] =
[z_{ni},x_{nk}] = [z_{nj},x_{nk}] = [z_{nj},x_{ni}]$.

Next, we need to show that $\bar{s_z}$ respects the defining
relations of $\Astar[u]$: \eqs{unit}{disjoint}, and the restricted
versions of \eqs{sony}{sons}. Relation \eq{unit} holds trivially.
For \eq{transitive} and \eq{transitivecomm}, choose $l\neq i,j,k$.
Then $\bar{s_z}$ acts on $x_{ij},x_{jk},x_{ik}$ as conjugation by
$z_{nl}$, so the equality is preserved. In order to prove that
$\bar{s_z}$ preserves \eq{disjoint} (for certain $i,j,k,l$), we
choose to write $\bar{s_z}(x_{ij}) = z_{nk} x_{ij} z_{nk}^{-1}$
and $\bar{s_z}(y_{kl}) = z_{ni} y_{kl} z_{ni}^{-1}$. Then we have
$[\bar{s_z}(x_{ij}),\bar{s_z}(y_{kl})] = [z_{nk} x_{ij}
z_{nk}^{-1}, z_{ni} y_{kl} z_{ni}^{-1}] = z_{nk} [x_{ij}, z_{ik}
y_{kl} z_{ki}] z_{nk}^{-1} = 1$ by \eq{xyzeq} (which holds under
our assumptions on $t,n$).
In order to check that
$$\bar{s_z}(s_x) \bar{s_z}(y_{ij}) \bar{s_z}(s_x)^{-1} = \bar{s_z}(x_{nk})\bar{s_z}(y_{ij})\bar{s_z}(x_{nk})^{-1}$$
(where $x < z$, $y \in X$ is arbitrary and $k \neq i,j$), we need
to choose parameters for the action of $\bar{s_z}$ on various
generators. We may assume that $i\neq j$, so for example $j \neq
n$. Then, write $\bar{s_z}(s_x) = [z_{nk},x_{nj}]s_x$,
$\bar{s_z}(y_{ij}) = z_{nk} y_{ij} z_{kn}$, and $\bar{s_z}(x_{nk})
= z_{nj} x_{nk} z_{jn}$. Then the equality becomes
$[z_{nk},x_{nj}] s_x z_{nk} y_{ij} z_{kn} s_x^{-1} [x_{nj},z_{nk}]
= z_{nj} x_{nk} z_{jk} y_{ij} z_{kj} x_{kn} z_{jn}$. We now choose
to write $s_x z_{nk}s_x^{-1} = x_{nj} z_{nk} x_{jn}$ and $s_x
y_{ij} s_x^{-1} = x_{nk} y_{ij} x_{kn}$, so the equation becomes
$$z_{jk} x_{nk} y_{ij} x_{kn} z_{kj} = x_{nk} z_{jk} y_{ij} z_{kj} x_{kn}.$$

There are two cases to consider. If $i = n$ then choose $r \neq
j,k,n$, and compute that
\begin{eqnarray*}
z_{jk} x_{nk} y_{ij} x_{kn} z_{kj}
& = & z_{jr} z_{rk} x_{nk} y_{ij} x_{kn} z_{kr} z_{rj} \\
& \explain{=}{\eq{xyzeqinn}} & z_{jr} x_{nk} y_{ij} x_{kn} z_{jr} \\
& = & x_{nk} z_{jr} y_{ij} z_{rj} x_{kn} \\
& = & x_{nk} z_{jr} z_{rk} y_{ij} z_{kr} z_{rj} x_{kn} \\
& = & x_{nk} z_{jk} y_{ij} z_{kj} x_{kn}.
\end{eqnarray*}
If $i \neq n$ then $x_{kn}$ commute with $y_{ij}$, and the
equality is equivalent to $z_{jk} [y_{ij},z_{kj} x_{nk} z_{jk}]
z_{kj} = 1$, which again holds by \eq{xyzeqinn}.

\forget
$$z_{kj} x_{nk} y_{ij} x_{kn} z_{jk} = x_{nk} z_{kj} y_{ij} z_{jk} x_{kn}.$$
There are two cases to consider. If $i = n$ then choose $r \neq
j,k,n$, and compute that
\begin{eqnarray*}
z_{kj} x_{nk} y_{ij} x_{kn} z_{jk}
& = & z_{rj} z_{kr} x_{nk} y_{ij} x_{kn} z_{rk} z_{jr} \\
& \explain{=}{\eq{xyzeq}} & z_{rj} x_{nk} y_{ij} x_{kn} z_{jr} \\
& = & x_{nk} z_{rj} y_{ij} z_{jr} x_{kn} \\
& = & x_{nk} z_{kj} z_{rk} y_{ij} z_{kr} z_{jk} x_{kn} \\
& = & x_{nk} z_{kj} y_{ij} z_{jk} x_{kn}.
\end{eqnarray*}
If $i \neq n$ then $x_{kn}$ commute with $y_{ij}$, and the
equality is equivalent to $z_{kj} [y_{ij},z_{jk} x_{nk} z_{kj}]
z_{jk} = 1$, which again holds by \eq{xyzeq}. \forgotten

Finally, we need to check \eq{sons}, where $x,y<z$, that is, to
show that $[\bar{s_z}(s_x),\bar{s_z}(s_y)] =
[\bar{s_z}(x_{ni}),\bar{s_z}(y_{nj})]$. Let $k \neq i,j,n$. We
choose to write $\bar{s_z}(s_x) = [z_{nk},x_{ni}]s_x$,
$\bar{s_z}(s_y) = [z_{nk},y_{nj}]s_y$, $\bar{s_z}(x_{ni}) = z_{nk}
x_{ni} z_{kn}$, and $\bar{s_z}(y_{nj}) = z_{nk} y_{nj} z_{kn}$, so
that $\bar{s_z}([x_{ni},y_{nj}]) = z_{nk} [x_{ni},y_{nj}] z_{kn}$.

Now the equation becomes
$$[z_{nk},x_{ni}]s_x  [z_{nk},y_{nj}]s_y s_x^{-1}[x_{ni},z_{nk}] s_y^{-1} [y_{nj},z_{nk}] = z_{nk} [x_{ni},y_{nj}] z_{kn}.$$
Letting $s_x$ act as conjugation by $x_{ni}$ on $z_{nk},y_{nj}$,
and $s_y$ as conjugation by $y_{nj}$ on $z_{nk},x_{ni}$, we get
$$[z_{nk},x_{ni}]x_{ni} [z_{nk},y_{nj}] x_{in} [s_x,s_y] y_{nj}[x_{ni},z_{nk}] y_{jn} [y_{nj},z_{nk}] = z_{nk} [x_{ni},y_{nj}] z_{kn},$$
which is equivalent (using only \eqs{unit}{transitivecomm}) to
$$y_{jn} x_{in} [s_x,s_y] y_{nj} x_{ni}   = 1.$$
This is relation \eq{sons}, so we are done.

Now that $\bar{s_z}$ is well defined, we are done by noting that
$\Astar[z]$ is generated from $\Astar[u]$ by adding one generator,
$s_z$, and the relevant portion of relations \eqs{sony}{sons},
which can be rewritten in the form $s_z w s_z^{-1} = \bar{s_z}(w)$
for $w = x_{ij}$ ($x \in X$) or $w = s_x$ ($x < z$). But these are
precisely the relations defining the semidirect product of $\Z =
\sg{s_z}$ acting on $\Astar[u]$ via the automorphism $\bar{s_z}$.
\end{proof}

\subsection{Action of $S_n$}\label{Snactss}

The symmetric group $S_n$ naturally acts on $\Fstar = (\F_t)^n$ by the action on indices:
\begin{equation}\label{actF}
\s^{-1} x_i \s = x_{\s(i)}.
\end{equation}
Likewise, it acts on
$A_{t,n}$ by
\begin{equation}\label{act}
\s^{-1} x_{ij} \s = x_{\s i,\s j}.
\end{equation}
(recall that $( \s \tau)(i) = \tau(\s(i))$).

\begin{rem}
Under the assumptions of Theorem \ref{AisomF}, the isomorphism
$F_{t,n} \isom A_{t,n}$ carries $x_{ij}$ to $x_{j}^{-1} x_{i} \in
F_{t,n}$, so it agrees with the actions of $S_n$ on both groups.
In particular the resulting semidirect products
$\semidirect{S_n}{F_{t,n}}$ and $\semidirect{S_n}{A_{t,n}}$ are
isomorphic.
\end{rem}

\subsection{Some identities in $A_{t,n}$}\label{identities}

The identities we prove here are most easily derived from the isomorphism
$A_{t,n} \isom F_{t,n}$, when it holds.
However, since we also want to cover the case $n=4$ and $t=3$,
our proofs use direct computation in $A_{t,n}$.

\begin{rem}\label{autoinv}
Choose any $x\in X$. The map fixing $y_{ij}$ for every $y \neq x$ and sending $x_{ij}\mapsto x_{ji}$
is an automorphism of $A_{t,n}$.
\end{rem}
\begin{proof}
This map preserves the relations \eq{unit} and \eq{disjoint}, and switches
\eq{transitive} and \eq{transitivecomm}.
\end{proof}

\begin{prop}\label{relcomm}
Suppose $n\geq 4$.
For every $x,y \in X$ and distinct $i,j,k$, the following relations hold in $A_{t,n}$.
\begin{eqnarray}
x_{jk}y_{ki}x_{ij} = y_{ji}x_{ik}y_{kj} \label{comm} \\
x_{jk}y_{ik}x_{ij} = y_{ij}x_{ik}y_{jk} \label{commi} \\
x_{kj}y_{ik}x_{ji} = y_{ij}x_{ki}y_{jk} \label{commii}
\end{eqnarray}
\end{prop}
\begin{proof}
Let $r \neq i,j,k$, and compute the quotient of the left hand and right hand sides in \eq{comm}:
\begin{eqnarray*}
x_{jk}y_{ki}x_{ij} (y_{ji}x_{ik}y_{kj})^{-1}
& = & x_{jk}y_{ki}x_{ij} y_{jk} x_{ki} y_{ij} \\
& = & x_{jk}y_{ki}x_{ir}x_{rj} y_{jk} x_{ki} y_{ij} \\
& = & x_{jk}x_{rj}y_{ki} y_{jk}x_{ir} x_{ki} y_{ij} \\
& = & x_{rk}y_{ji} x_{kr} y_{ij} = 1.
\end{eqnarray*}

The other two relations are obtained by inverting $x_{rs}$, or $x_{rs}$ and $y_{rs}$,
using Remark \ref{autoinv}.
\end{proof}

\begin{prop}\label{xyz2}
Suppose $n \geq 5$ or $t \leq 2$.

If $u,v,w\in X$ and $i,j,k,s$ are distinct, we have that
\begin{equation}\label{xyzeq2}
u_{si}v_{ij}u_{js} w_{sk} = w_{sk} u_{ki} v_{ij} u_{jk}.
\end{equation}
\end{prop}
\begin{proof}
This identity is equivalent to \Eq{xyzeq}, as the following computation
shows.
\begin{eqnarray*}
u_{si}v_{ij}u_{js} w_{sk} u_{kj}v_{ji}u_{ik} w_{ks}
    & = & u_{si}v_{ij}u_{js} w_{sk} (u_{kj}v_{ji}u_{jk})u_{ij} w_{ks} \\
    & \explain{=}{\eq{xyzeq}} & u_{si}v_{ij}u_{js} u_{kj}v_{ji}u_{jk} w_{sk} u_{ij} w_{ks} \\
    & = & u_{si}v_{ij} u_{ks} v_{ji} u_{jk} w_{sk} u_{ij} w_{ks} \\
    & = & u_{si} u_{ks}  u_{jk} u_{ij} = 1.
\end{eqnarray*}
\end{proof}

\section{The group $\CoxY{T}$}\label{mainsec}

\subsection{The main result}

\forget Let $M = A_{1,n}$ be the Abelian group generated by
$x_{ij}$ ($1 \leq i, j \leq n$), subject to the relations $x_{ii}
= 1$ and $x_{ij}x_{jk} = x_{ik}$. The symmetric group acts on $M$
by $\s^{-1} x_{ij} \s = x_{\s(i)\s(j)}$. A key observation in the
identification of $\CoxY{T}$ is the following. Fix an edge $x =
(a,b) \in T$. We can alter the map $\phi \co \CoxY{T}\ra S_n$
defined in Section \ref{Snpres} by setting $x \mapsto (ab)x_{ab}$.
This map is still well defined from $\CoxY{T}$ to
$\semidirect{S_n}{M}$. Moreover, if $x = u_1$ belongs to a cycle
$u_1,u_2,\dots,u_m \in T$, then $(u_1 \dots u_{m-1})(u_2\dots
u_m)^{-1} \mapsto x_{ab}$, so the map is onto. \forgotten

In this subsection we prove the following
\begin{thm}\label{main}
Let $T$ be a connected graph on $n$ vertices.
Then $\CoxY{T}$ is isomorphic to
the semidirect product $G = \semidirect{S_n}{A_{t,n}}$, where $t$ is the rank of the free group
$\pi_1(T)$.
\end{thm}

Recall that $\CoxY{T}$ was defined in Definition \ref{CYT}, the
group $A_{t,n}$ in Subsection \ref{Atnss}, and the action of $S_n$ on $A_{t,n}$ in Subsection \ref{Snactss}.

The proof is rather direct: we define maps  $\Phi \co \CoxY{T} \ra
\semidirect{S_n}{A_{t,n}}$ and  $\Psi \co
\semidirect{S_n}{A_{t,n}} \ra \CoxY{T}$, show that they are well
defined, and check that they invert each other. Showing that
$\Psi$ is well defined is the longest part of the proof, relying
on Corollary \ref{PsiPreserveAct}.

Choose a spanning tree $T_0 \sub T$,
and let $X = T-T_0$.
Arbitrarily fix a direction for every edge $x \in X$, so $x$ has a starting
point
$a$ and an ending point $b$.
There is a unique path in $T_0$ from $a$ to $b$, and
together with $x$ this path forms a cycle in $T$, called a \defin{basic cycle}
(this is the unique cycle in $T_0 \cup \set{x}$).
The set of basic cycles is a free set of generators for the fundamental
group $\pi_1(T)$.  Let $t = |X|$ be the number of basic cycles.
We will always use $X$ as the set of size $t$ indexing the generators of
$A_{t,n}$.

Note that the number $t$ of basic cycles in a graph on $n = 4$ vertices is
bounded by $t \leq 3$, and if $n = 3$ there is
at most one cycle. If $t = 0$ Theorem \ref{main} is vacuous
(since $A_{0,n} = 1$ and $\CoxY{T} = S_n$ by Corollary \ref{Artree}),
so we may assume $t \geq 1$ and $n \geq 3$.
It follows that the only case where Theorem \ref{AisomF} does not apply is
where $n = 4$ and $t = 3$, which is the case iff
$T$ is the complete graph $K_4$ on four vertices (this is the only graph on
$4$ vertices with $3$ independent cycles).

\begin{defn}[Definition of $\Phi: \CoxY{T} \ra \semidirect{S_n}{A_{t,n}}$]\label{Phidef}
Let $u \in T$ be an edge from $a$ to $b$. We set
$$\Phi(u) = \begin{cases}(ab) & \mbox{if $ u \in T_0$,} \\ (ab)u_{ab} & \mbox{if $u \in X$.} \end{cases}$$
\end{defn}

\begin{prop}
The map $\Phi$ is well defined on $\CoxY{T}$.
\end{prop}
\begin{proof}
We need to check the relations \eqs{square}{fork}.
If $u \in T_0$ then $\Phi(u)^2 =
(ab)^2 = 1$, and for $u \in X$ we have $\Phi(u) =
(ab)u_{ab}(ab)u_{ab} = (ab)^2 u_{ba}u_{ab} = u_{aa} = 1$ by \eq{unit},
so \Eq{square} is preserved.

Let $u,v$ be two disjoint edges. If both belong to $T_0$ then obviously
$\Phi(u)$ and $\Phi(v)$ commute. If $u=(a,b) \in T_0$ and $v = (c,d) \in X$,
then $\Phi(u) = (ab)$ commutes with $\Phi(v) = (cd)x_{cd}$ since $(ab)$
acts trivially on $x_{cd}$. And if $u,v \in X$, $(ab)u_{ab}$ commutes
with $(cd)v_{cd}$ by \eq{disjoint}. Thus \eq{commute} is preserved.

Let $u,v$ be two intersecting edges, \eg\ $u = (a,b)$ and $v = (b,c)$.
If $u,v \in T_0$ then \Eq{triple}
is satisfied by the transpositions $(ab),(bc)$. If $u \in T_0$
and $v \in X$ then $\Phi(u)\Phi(v) = (ab)(bc)v_{bc} = (cba)v_{bc}$ and
\begin{eqnarray*}
(cba)v_{bc}(cba)v_{bc}(cba)v_{bc} & = & (cba)^2 v_{ab}v_{bc}(cba)v_{bc} \\
    & = & (cba)^2v_{ac}(cba)v_{bc} = v_{cb}v_{bc} = 1.
\end{eqnarray*}
 The last case to check is when $u,v \in X$ (this can only happen if $t \geq 2$, so we assume $n\geq 4$).
Then
$\Phi(u) = (ab)u_{ab}$ and $\Phi(v) = (bc)v_{bc}$
so we have $\Phi(u)\Phi(v)\Phi(u) =
(ab)u_{ab}(bc)v_{bc}(ab)u_{ab} = (ac)u_{bc}v_{ac}u_{ab}$ and
$\Phi(v)\Phi(u)\Phi(v) = (bc)v_{bc}(ab)u_{ab}(bc)v_{bc} =
(ac)v_{ab}u_{ac}v_{bc}$, and the equality follows from \eq{commi}.
Here we assumed that the ending point
of $u$ is the starting point of $v$.
The other possible ways to direct $u,v$ are checked similarly
using \eq{comm} and \eq{commii}, so \Eq{triple} is verified.

Finally, we need to check that $\Phi$ respects \eq{fork}.
Suppose $u,v,w \in T$ meet in the
same vertex $s$, which $u$ connects to $i$, $v$ to $j$ and $w$
to $k$.
If $u,v \in T_0$ then we are done since $(is)(js)(is) =
(ji)$ commutes with $(ks)$ and $(ks)w_{ks}$.
If $u \in T_0$ and $v \in X$ then $\Phi(u)\Phi(v)\Phi(u) =
(is)(js)v_{js}(is) = (ji)v_{ji}$ which commutes with $(ks)$ and
$(ks)w_{ks}$ by Relations \eq{disjoint} and \eq{act}. The last
case to consider is when $u,v,w\in X$.
This cannot happen if $n \leq 3$, or if $n=4$ and $t = 3$, for in this case $T = K_4$ and
the complement of $u,v,w$ would be a triangle, not containing a spanning subtree.
So we may assume that $n\geq 5$, and then Proposition \ref{xyzredun} applies.
Then $\Phi(u)\Phi(v)\Phi(u) =
(ij)u_{sj}v_{ji}u_{is}$, and
\begin{eqnarray*}
\, [\Phi(u)\Phi(v)\Phi(u),\Phi(w)]
    & = & (ij)u_{sj}v_{ji}u_{is} (ks) w_{ks} (ij) u_{sj}v_{ji}u_{is} (ks) w_{ks} \\
    & = & u_{si}v_{ij}u_{js} w_{sk} u_{kj}v_{ji}u_{ik} w_{ks} = 1
\end{eqnarray*}
by \Eq{xyzeq2}.
\end{proof}

\begin{rem}\label{restPhi}
By definition, the restriction of $\Phi \co \CoxY{T} \ra
\semidirect{S_n}{A_{t,n}}$ to the subgroup $\sg{T_0}$ is $\Phi(u)
= (ij)$ (where $i,j$ are the vertices of $u$), so $\Phi$ is an
extension  of the isomorphism $\phi_0 \co \sg{T_0} \ra S_n$ of
Proposition \ref{parabolicSn}.
\end{rem}

\begin{defn}[Definition of $\Psi: \semidirect{S_n}{A_{t,n}} \ra \CoxY{T}$]
We define $\Psi$ on $S_n$ to invert $\Phi$:
$\Psi(\s)$ is the unique element of $\sg{T_0}$ which $\Phi$ takes to $\s$.

If $x \in X$, let $u_2,\dots,u_m \in T_0$ be the path connecting
the vertices of $x$; label the vertices on the cycle so that $u_i$
connects $i-1,i$ (and $x$ connects $1,m$), and define $\gamma_i$,
$i = 1,\dots,n$, as in \eq{gammadef} (extending the definition
\eq{gammacycdef}). Then define (for $i,j = 1,\dots,n$)
$$\Psi(x_{ij}) = \gamma_j^{-1}\gamma_{i}.$$
\end{defn}

\begin{prop}
The map $\Psi$ is well defined on $\semidirect{S_n}{A_{t,n}}$.
\end{prop}
\begin{proof}
As $\Psi$ is well defined on $S_n$  (see Remark \ref{restPhi}), we
only need to check the defining relations: \eqs{unit}{disjoint}
and \eq{act}. Observe that \eq{unit} and \eq{transitivecomm} are
trivially preserved. Next, we already proved \eq{act} in Corollary
\ref{PsiPreserveAct}. This equation is very useful, as it allows
us to carry 'local' proofs (for indices of our choice) to the
general case.

In order to check Equation \eq{transitive},
we first prove that $\Psi(x_{1m}) = \gamma_m^{-1} \gamma_1$ and $\Psi(x_{2m}) = \gamma_m^{-1} \gamma_2$ commute.
Recall from \Eq{gg} that if $i,j$ are on the cycle, we have
$\gamma_i\gamma_j^{-1} = \gamma_{j+1}^{-1}\gamma_{i+1}$. Now compute:
\begin{eqnarray*}
\, [\gamma_m^{-1} \gamma_1, \gamma_m^{-1} \gamma_2]
    & = & \gamma_m^{-1} \gamma_1 \gamma_m^{-1} \gamma_2 \gamma_1^{-1} \gamma_m \gamma_2^{-1} \gamma_m \\
    & = & \gamma_m^{-1} \gamma_1^{-1} \gamma_2 \gamma_2 \gamma_1^{-1} \gamma_3^{-1} \gamma_1 \gamma_m \\
    & = & \gamma_m^{-1} \gamma_1^{-1}  \gamma_2 \gamma_2^{-1} \gamma_3 \gamma_3^{-1}  \gamma_1 \gamma_m = 1.
\end{eqnarray*}
Now, given any distinct $i,j,k$, choose a permutation $\s$ which carries $1,m,2$ to $i,j,k$;
then using \eq{transitivecomm} we have
\begin{eqnarray*}
\Psi(x_{ij})\Psi(x_{jk}) & = & \s^{-1} \Psi(x_{1m}) \Psi(x_{m2}) \s \\
    & = & \s^{-1} \Psi(x_{m2}) \Psi(x_{1m}) \s = \s^{-1} \Psi(x_{12}) \s = \Psi(x_{ik}).
\end{eqnarray*}

If $t = 1$ then \Eq{disjoint} is trivial, so we are done in this case.
Using Proposition \ref{PsiPhi} and Corollary \ref{PhiPsi}
(which are independent of the current proposition),
the proof of Theorem \ref{main} itself is complete if $t = 1$.
This will be used instead of a
lengthy case analysis below.

We will now show that $\Psi$ respects Equation \eq{disjoint}. Let
$x\neq y \in X$. Let $\gamma_i$ be defined as usual with respect
to $x$, and similarly define $\gamma'_{i'}$ for $y$ (we label the
vertices of the basic cycle of $x$ as $1,\dots,m$ for $x$, and the
vertices of the basic cycle of $y$ as $1',\dots,m'$). Let $g_x =
\gamma_m^{-1}\gamma_1 x$ and $g_y =
{\gamma'}_{m'}^{-1}\gamma'_{1'} y$; by definition of the
$\gamma_i, {\gamma'}_{i'}$, we have that $g_x,g_y \in \sg{T_0}$.
First assume $x,y$ are disjoint, so the endpoints $1,m$ of $x$ and
$1',m'$ of $y$ are distinct. Let $T_1 = T_0 \cup \set{x}$. Working
in the abstract group $\CoxY{T_1}$, we see that $x$ commutes with
$g_y$ (since their images under the appropriate $\Psi$ are $(1
m)x_{1m}$ and $(1'm')$). By Remark \ref{parabolicweak}, $x$
commutes with $g_y$ in $\sg{T_0,x} \leq \CoxY{T}$. Similarly, $y$
commutes with $g_x$. Moreover, $g_x$ commutes with $g_y$, as they
correspond to disjoint transpositions in $\sg{T_0} = S_n$, and
finally $x,y$ commute by assumption. Since $g_x,x$ commute with
$g_y,y$, we have that $\Psi(x_{1m}) = \gamma_m^{-1} \gamma_1 = g_x
x$ commutes with $\Psi(y_{1'm'}) = g_y y$. Conjugating by
arbitrary $\s \in S_n$, we have that $\Psi(x_{ij})$ commutes with
$\Psi(y_{kl})$ for arbitrary disjoint $i,j,k,l$, which is what we
need.

Now assume that $x,y$ touch each other; let $1,m$ denote the
endpoints of $x$, and $1,m'$ the endpoints of $y$. Consider the
paths in $T_0$ connecting $1$ to $m$ and $m'$. Since $T_0$ is a
tree, their intersection is connected; we denote by
$u_2,\dots,u_m$ the path from $1$ to $m$, and by
$v_2,\dots,v_{m'}$ the path from $1$ to $m'$ (the two lists
coincide at the beginning). See Figure \ref{pic2}. As always,
denote by $m-1$ the other vertex of $u_m$. Exchanging $x,y$ if
necessary, we may assume that $u_m$ is not in
$\set{v_2,\dots,v_{m'}}$. In particular $m-1 \neq m'$.

\begin{figure}[!h]
\ifXY
\begin{equation}\nonumber
\xymatrix{ 
{}
    & m \ar@{=}@/^2pt/[r]^{x} \ar@{-}@/_4pt/[ld]^{u_m}
    & 1 \ar@{=}@/^2pt/[r]^{y} \ar@{-}[d]^{v_2}_{u_2}
    & m' \ar@{.}@(r,d)@/^14pt/ [ddr]
    &   \\
m-1 \ar@{.}[d]
    &
    & 2 \ar@{.}[d]
    &
    &   \\
\ar@{.}@(d,u)@/_42pt/ [rr]
    &
    & \ar@{.}@(d,u)@/_36pt/ [rr]
    &
    &   \\
{}
    &
    &
    &
    &   \\
}
\end{equation}
\else

\bigskip

(A picture of the circle).

\bigskip
\fi 
\caption{}\label{pic2}
\end{figure}

Let $u = u_2 \dots u_{m-1}$, and $T_2 = \set{u_2,\dots,u_m,x,y}
\sub T$. This subgraph has one cycle, and in $\CoxY{T_2}$ (which
is isomorphic to $\semidirect{S_{m+1}}{\Z^{m}}$ by the case $t =
1$ of Theorem \ref{main}), one checks that $y$ commutes with
$u^{-1} g_x x u$, as they correspond to $(1m')$ and $x_{m\,m-1}$
(note that in this subgraph, $y$ is part of the spanning subtree).
Similarly let $T_3 = \set{u_2,\dots,u_m,x,v_2,\dots,v_{m'}}$; this
subgraph too has one cycle, and in $\CoxY{T_3}$ we have that $g_y
= v_2 \dots v_{m'} \dots v_2$ commutes with $u^{-1} g_x x u$ (as
they again correspond to $(1m')$ and $x_{m\,m-1}$). Finally
$u^{-1} g_x x u = u^{-1} \gamma_m^{-1} \gamma_1 u =
\gamma_{u(m)}^{-1}\gamma_{u(1)} = \gamma_m^{-1} \gamma_{m-1} =
\Psi(x_{m-1\,m})$ by \Eq{doubleact}, so $\Psi(x_{m-1\,m})$
commutes with $\Psi(y_{1m'}) = \gamma_{m'}^{-1} \gamma_{1} = g_y
y$. Conjugating by arbitrary $\s \in S_n$ we see that $\Psi$
respects \eq{disjoint}.
\end{proof}

\begin{prop}\label{PsiPhi}
The composition $\Psi \circ \Phi$ is the identity on $\CoxY{T}$.
\end{prop}
\begin{proof}
Let $u=(i,j) \in T_0$. Then $\Phi(u) = (ij)$, and by definition $\Psi \Phi(u) =
\Psi(ij) = u$. Now let $x \in X$, and set $u_1 = x$.
As usual we label the vertices along the basic cycle attached to $x$,
so that $x$ connects $1$ and $m$, and by definition $\Phi(x) = (1m)x_{1m}$.
Let
$\gamma_i$ be defined as in \eq{gammadef},
and compute that $\gamma_m^{-1} \gamma_1 u_1 = u_2 \dots u_m \dots u_2$
is an element of $\sg{T_0}$ which maps to $(1m)$, so $\Psi((1m)) = u_1 \gamma_1^{-1} \gamma_m$.
By definition $\Psi(x_{1m}) = \gamma_m^{-1}\gamma_1$, so that $\Psi\Phi(x) = \Psi((1m)x_{1m}) = u_1 = x$,
as required.
\end{proof}

For the other direction, we first compute $\Phi(\gamma_a)$.
\begin{prop}\label{Phigamma} 
Let $\gamma_a$ be defined as in \eq{gammadef},
and set $\tau = (m \dots 3 2 1)$.
Then for every $a = 1, \dots ,n$ we have that $\Phi(\gamma_a) = \tau x_{am}$.
\end{prop}
\begin{proof}
First compute $\Phi(\gamma_j)$ where $j$ is on the cycle (that is, $1 \leq j \leq m$).
We may use \eq{gammacycdef} as the definition of $\gamma_j$. If $j < m$
we have
\begin{eqnarray*}
\Phi(\gamma_j) & = & \Phi(u_{j+2} \dots u_m u_1 \dots u_j) \\
    & = & (j+1\,j+2)\dots (m-1\,m)(1m)x_{1m}(12)\dots(j-1\,j) = \tau x_{jm},
\end{eqnarray*}
and for $j = m$, $\Phi(\gamma_m) =\Phi(u_{2} \dots u_m) = (12)\dots(m-1\,m)=\tau = \tau x_{mm}$.

In general, $\gamma_a$ is defined in \eq{gammadef}.
Let $v \in T_0$ be an edge connecting a vertex $j$ on the cycle to a vertex $b$ (not on the cycle).
By definition, $\til{v} = u_{j+1} v u_{j+1}$. So if
$j < m$ we have $\Phi(\til{v}) = (j\,j+1)(ib)(j\,j+1) = (j+1\,b)$,
while if $j = m$ we have $\Phi(\til{v}) = (1m)x_{1m}(mb)(1m)x_{1m} = (1b)x_{1b}$.

Now let $a$ be an arbitrary vertex, and let $w_1,\dots,w_s$ denote the path connecting $a$ to a vertex $j$ on the cycle.
Write $w = w_2 \dots w_s$, and let $b = w_1(j)$. Note that $w(b) = a$.
Since $w$ does not touch the cycle we have that $\til{w} = w$, so
by definition $\gamma_a = w^{-1} \til{w_1} \gamma_j w_1 w$. Now if $j<m$,
\begin{eqnarray*}
\Phi(\gamma_a)
    & = & \Phi(w^{-1} \til{w_1} \gamma_j w_1 w) \\
    & = & \Phi(w)^{-1} \Phi(\til{w_1} \gamma_j w_1) \Phi(w) \\
    & = & \Phi(w)^{-1} (j+1\,b) \tau x_{jm} (j b) \Phi(w) \\
    & = & \Phi(w)^{-1} (j+1\,b) \tau (j b) \Phi(w) x_{am} \\
    & = & \Phi(w)^{-1} \tau \Phi(w) x_{am} = \tau x_{am},
\end{eqnarray*}
where the last equality follows since $w$ does not intersect the cycle.

Similarly if $j = m$, we have
\begin{eqnarray*}
\Phi(\gamma_a)
    & = & \Phi(w^{-1} \til{w_1} \gamma_m w_1 w) \\
    & = & \Phi(w)^{-1} \Phi(\til{w_1} \gamma_m w_1) \Phi(w) \\
    & = & \Phi(w)^{-1} (1 b) x_{1b} \tau (m b) \Phi(w) \\
    & = & \Phi(w)^{-1} (1 b) \tau x_{m b} (m b) \Phi(w) \\
    & = & \Phi(w)^{-1} (1 b) \tau (m b) \Phi(w) x_{a m} \\
    & = & \Phi(w)^{-1} \tau \Phi(w) x_{a m} = \tau x_{a m}.
\end{eqnarray*}
\end{proof}

\begin{cor}\label{PhiPsi}
$\Phi\circ \Psi$ is the identity on $\semidirect{S_n}{A_{t,n}}$.
\end{cor}
\begin{proof}
For $\s \in S_n$ we have that $\Psi(\s) \in \sg{T_0}$ and by definition $\Psi\Phi(\s) = \s$.
Let $x \in X$ and let $i,j = 1,\dots,n$. Then by definition $\Psi(x_{ij}) = \gamma_j^{-1} \gamma_i$ where
$\gamma_i,\gamma_j$ are constructed as usual.

By the last proposition and \eq{transitivecomm},
$\Phi(\gamma_j^{-1}\gamma_i) = x_{mj}\tau^{-1} \cdot \tau x_{im} = x_{ij}$.
\end{proof}

This ends the proof of Theorem \ref{main}.

\subsection{The map $u \mapsto \til{u}$}\label{tildess}

As always, let $T_0$ be a spanning subtree of $T$ and $x \in T-T_0$, and let $u_2,\dots,u_m \in T_0$ be the completion
of $x$ to a basic cycle in $T$.
In Section \ref{actionsec} we defined a map $\tilmap$ ($u \in T$), and showed
it induces a homomorphism from $\sg{T_0}$ to $\CoxY{T}$.
This map has an interesting role in the proof of Theorem \ref{main},
so before moving to the applications, we will prove the following:
\begin{prop}
Assume $T \neq K_4$.
Then the map $\tilmap$ of Definition \ref{tildedef} extends to an automorphism of $\CoxY{T}$.
\end{prop}
\begin{proof}
We identify the groups $\Astar$ and $\Fstar$ (using the proof of
Theorem \ref{AisomF}). From Theorem \ref{main} we have the
isomorphism $\Psi \co \semidirect{S_n}{F_{t,n}} \ra \CoxY{T}$,
inverted by $\Phi$ of Definition \ref{Phidef}. Set $\tau = (m
\dots 3 2 1)$. We claim that for every $u \in T$,
$$\til{u} = \Psi(\tau x_m^{-1} \cdot \Phi(u) \cdot x_m \tau^{-1}).$$
The right hand side is an inner automorphism of $\semidirect{S_n}{\Fstar}$, induced by the element $\tau x_m^{-1}$.
In order to prove this claim, we choose for every $u \in T$ a vertex $c$ not lying on $u$. Then $x_c$ commutes with $\Phi(u)$, and
we need to prove
$$\til{u} = \Psi(\tau x_m^{-1} x_c \cdot \Phi(u) \cdot x_c^{-1} x_m \tau^{-1}),$$
which using $\Psi(x_m^{-1} x_c) = \gamma_m^{-1} \gamma_c$
and the fact that $\Psi(\tau) = \gamma_m$ (see Proposition \ref{Phigamma}),
translates to
$$\til{u} = \gamma_c u \gamma_c^{-1}.$$
Note that apart from the assumption that $u_2,\dots,u_m \in T_0$, the map $\tilmap$ is
independent of $T_0$. This is not true for the $\gamma_i$ in general (see \Eq{gammadef}), but for $i = 1,\dots,m$, we
do have that $\gamma_i$ is independent of $T_0$ (see \eq{gammacycdef}). Now, for every $u \in T$ which
does not touch the cycle twice, we can a priori choose $T_0$ to include $u$; then the result follows from Proposition \ref{singleactprop}.
The remaining case is where $u \in T$ touches the cycle twice, but we have shown in Case \eq{touchtwice}
of the calculations preceding Definition \ref{tildedef} that $\til{u} = \gamma_{i-1} u \gamma_{i-1}^{-1}$ in this case.
\end{proof}

In contrast, let us consider the case $T = K_4$. We label the vertices and
edges as in Figure \ref{K4gr}.

\begin{figure}[!h]
\ifXY
\begin{equation}\nonumber
\xymatrix@C=30pt@R=18pt{
{}
    & 4 \ar@{-}[ddr]|{\,z\,} \ar@{-}[ddl]|{\,y\,} \ar@{-}[d]|{\,v\,}
    &
\\
{}
    & 2 \ar@{-}[dr]|{\,u_3\,} \ar@{-}[dl]|{\,u_2\,}
    &
\\
1 \ar@{-}[rr]|{\,x\,}
    &
    & 3
}
\end{equation}
\else
\bigskip

(The $K_4$ graph).

\bigskip
\fi 
\caption{$T = K_4$}\label{K4gr}
\end{figure}

We choose the spanning subtree $T_0 = \set{u_2,u_3,v}$, and
consider the basic cycle of $x$ (oriented as $u_1 = x, u_2, u_3$).
Then by definition of $\tilmap$ we have $\til{u_2} = u_3$,
$\til{u_3} = x$, $\til{x} = u_2$, $\til{v} = u_3 v u_3$, $\til{y}
= u_2 y u_2$, and $\til{z} = xzx$. It is easy to directly check
that $\sg{u_3, x, u_3 v u_3} \isom S_4$ (as Corollary \ref{ishom}
predicts), but we do not know if $[\til{v},\til{y}\til{z}\til{y}]
= [u_3 v u_3, u_2 y u_2xzx u_2 y u_2]$ is the trivial element in
$\CoxY{T}$. It is interesting to note that under the
identification of $\CoxY{T}$ with $\semidirect{S_4}{A_{3,4}}$,
this element evaluates to
$y_{32}(x_{41}z_{12}x_{24}y_{43}x_{32}z_{21}x_{13}y_{34})y_{23}$,
which is trivial iff \eq{xyzeq2} holds.

\section{Applications}\label{appl}

\subsection{The structure of $\CoxY{T}$}

In this section we apply Theorem \ref{main} to answer some natural
questions on the structure of $\CoxY{T}$ and its normal subgroup
$K(T) = \Ker(\maptoSn \co \CoxY{T} \ra S_n)$.

Let $t$ denote the rank of $\pi_1(T)$; $t = 0$ if $T$ is a tree, and $t = 1$ if $T$ has a unique cycle.
Recall that $K_4$ is the complete graph on $4$ vertices, the only exception to Theorem \ref{AisomF}.

\begin{cor}\label{whensolve}
Let $T$ be a connected graph on $n$ vertices.

a. $\CoxY{T}$ is virtually solvable iff $t\leq 1$.

b. If $T \neq K_4$, then $K(T) \sub \pi_1(T)^{n} = (\F_t)^{n}$ .

c. $K(T)$ maps onto $\pi_1(T)^{n-1}$.

d. The group $\CoxY{T}$ contains a subgroup isomorphic to $\pi_1(T)$.
\end{cor}
\begin{proof}
By Theorem \ref{main}, $\CoxY{T} \isom \semidirect{S_n}{A_{t,n}}$.
If $t = 0$ then the kernel $K(T) = \Ker(\maptoSn \co \CoxY{T} \ra
S_n)$ is trivial  (Corollary \ref{Artree}). If $t = 1$ we have
that $\CoxY{T} = \semidirect{S_n}{\Z^{n-1}}$ and $K(T) = \Z^{n-1}$
(Example \ref{A1n}). Here $\Z^{n-1}$ is the large irreducible
component of the standard representation of $S_n$. If $t \geq 2$
then $\CoxY{T}$ cannot be even virtually super-solvable, by Part
c.

Part b. follows from Theorems \ref{main} and \ref{AisomF}.
Recall the presentation of $\Fstar = \F_t^{n}$
from Subsection \ref{Ftnss}, and let $F_0$ denote the subgroup
generated by the elements $x_{i}$ for $i = 1,\dots,n-1$
(which is isomorphic to $(\F_t)^{n-1}$).
The map defined by $x_i \mapsto x_i$ for $i< n$ and $x_n \mapsto 1$ (for all $x \in X$) is onto,
since $x_i x_n^{-1} \in F_{t,n}$ covers $x_i$.
If $T \neq K_4$ we are done (as $K(T) \isom F_{t,n}$),
but even if $T = K_4$, $F_{t,n}$ is a quotient of
$A_{t,n} \isom K(T)$ by Theorem \ref{AisomF}.

Finally,
$x \mapsto x_1x_2^{-1}$ defines a map $\pi_1(T) \ra F_{t,n}$ which is obviously injective, thus
proving Part d. for $T \neq K_4$ (since then $F_{t,n} \sub \CoxY{T}$).
For $T = K_4$, $\F_3 \sub \F_2 \sub F_{2,4} \isom A_{2,4} \sub A_{3,4}$ by Proposition \ref{At1At},
so that the same result holds.
\end{proof}

Note that Part d. above can be strengthened to state that $\CoxY{T}$ contains $[\frac{n-1}{2}]$ commuting copies
of $\pi_1(T)$, using the free subgroups $\sg{x_{12} \suchthat x \in X}$, $\sg{x_{34} \suchthat x \in X}$,... .

\begin{exmpl}
Consider the graph $T$ of Figure \ref{sixpts}, where the spanning subgraph is $T_0 = \set{a,b,c,d,e}$ and $X = \set{x,y,z}$.
\begin{figure}[!h]
\ifXY
\begin{equation}\nonumber
\xymatrix@C=24pt@R=24pt{
1 \ar@{-}[r]^{a} \ar@{=}[d]^{x} \ar@{-}[rd]^{c}
    & 2 \ar@{-}[r]^{b} \ar@{-}[rd]^{d}
    & 3 \ar@{=}[d]^{y} \\
4 \ar@{-}[r]^{e}
    & 5 \ar@{=}[r]^{z}
    & 6 \\
}
\end{equation}
\else
\bigskip

(A graph with square and two triangles).

\bigskip
\fi 
\caption{}\label{sixpts}
\end{figure}

The isomorphism $K(T) \isom A_{3,6}$
gives $x_{14} = cecx$, $y_{36} = bdby$ and $z_{56} = cadacz$.
As three commuting free subgroups we can take $$\sg{caexecac,dadbybad,adczcd},$$
$$\sg{bacecxab,ecabdydace,bdecadaczedb}$$
and
$$\sg{daexecad,cadbybac,cadacz}.$$

\end{exmpl}

\forget 
\begin{cor}
If $T$ has at least two cycles, then the kernel of $\maptoSn \co
\Cox{T} \ra S_n$ maps onto $\F_t^{n-1}$.

In particular, $\Cox{T}$ cannot be virtually solvable.
\end{cor}
\begin{proof}
Immediate from Part b. of the last corollary, as $\Ker(\CoxY{T}\ra S_n)$ is a quotient of $\Ker (\Cox{T}\ra S_n)$.
\end{proof}
\forgotten

The following interesting property is immediate from Theorem \ref{main}:
\begin{rem}
The group $\CoxY{T}$ depends on the graph $T$ only through the
number of basic cycles and the number of vertices in the graph.
\end{rem}

Let $T$ be a connected graph. The Abelianization of $\Cox{T}$ is
easily seen to be $\Z/2\Z$ as every relation of the form $uvu =
vuv$ becomes $u = v$. The isomorphism $\CoxY{T} \isom
\semidirect{S_n}{A_{t,n}}$ proves that the commutator subgroup of
$\CoxY{T}$ coincides with $\semidirect{A_n}{A_{t,n}}$ (where $A_n$
is the alternating group). On the other hand we have
\begin{prop}
Let $T$ be a connected graph. The Abelianization of $K(T)$ is
$\Z^{t(n-1)}$.
\end{prop}
\begin{proof}
As an Abelian group, $A_{t,n}$ is freely generated by the $x_{i,i+1}$ ($x \in X$, $i = 1,\dots,n-1$).
\end{proof}

Another easy property of $F_{t,n}$ allows us to conclude the
following
\begin{prop}
Let $T \neq K_4$ be a connected graph. Then the kernel $K(T)$ of
$\maptoSn \co \CoxY{T} \ra S_n$ is torsion free.
\end{prop}

The following is also of some interest:
\begin{prop}
If $T \neq K_4$ is a connected graph, then $\CoxY{T}$ is residually finite.
\end{prop}
\begin{proof}
Given $1 \neq w \in \CoxY{T}$, we need to show that $w$ is outside a finite index normal subgroup of $\CoxY{T}$.
If $w \not \in K(T)$, we are done; and if $w \in K(T) \isom F_{t,n}$ use the residually finiteness of the free group.
\end{proof}

Since $\semidirect{S_n}{F_{t,n}} \sub \semidirect{S_n}{\Fstar}$,
we also have
\begin{prop}
The word problem is solvable in $\CoxY{T}$ for every connected $T
\neq K_4$.
\end{prop}

\subsection{Parabolic subgroups}\label{parabolicss}

In Proposition \ref{parabolicSn} we saw that a subgroup of $\CoxY{T}$ generated by a subtree $T_0$ is
isomorphic to the abstract group $\CoxY{T_0}$ defined on that tree. This can be generalized. We will show that
parabolic subgroups are well behaved,
even without assuming the smaller defining graph is connected.
\begin{prop}\label{parabolic}
Let $T$ be a graph, and $T' \sub T$ a subgraph.
Let $C = \sg{T'}$ be the subgroup of $\CoxY{T}$, generated by the vertices $u \in T'$.

If $T' \neq K_4$, then $C \isom \CoxY{T'}$. 
\end{prop}

This result should be compared to the situation for Coxeter groups.
If $D = \sg{s_1,\dots,s_k}$ is a Coxeter group, then it is well known (\eg\ \cite[Section
5.5]{Humph}) that for $I \sub \set{1,\dots,k}$, the subgroup of $D$ generated
by $\set{s_i \suchthat i \in I}$ is isomorphic to the appropriate
abstract Coxeter group.

Our result that  $\CoxY{T'}$ naturally embeds into $\CoxY{T}$ (as
the subgroup generated by the edges of $T'$) is of the same
nature; even better, if $T'$ is connected we show that $\CoxY{T'}$
is a retract of $\CoxY{T}$.

\begin{proof}[Proof of Proposition \ref{parabolic}]
First assume that $T'$ is connected. By induction we may assume
that $|T-T'| = 1$. Let $w$ denote the element in $T-T'$. Let $T_0
\sub T'$ be a spanning subtree (of $T$), $X' = T' - T_0$ and $X =
T-T_0$. Inspecting the image of $\sg{X'}$ under the isomorphism
$\Phi \co \CoxY{T} \llra \semidirect{S_n}{A_{X,n}}$, we only need
to show that the subgroup $\sg{x_{ij} \suchthat x\in X'}$ of
$A_{X,n}$ is isomorphic to $A_{X',n}$, and this is Proposition
\ref{At1At}.

Now assume that $T'$ is not connected. Again we apply induction but this time it is necessary
to choose the chain $T' = T^{0} \sub T^{1} \sub \dots \sub T^{s} = T$
of subgraphs (with $|T^{i+1}-T^{i}|=1$), in a way that $T^i \neq K_4$.
Indeed, assume $T^i = K_4$ for some $i$. Since $T' \neq K_4$, $i>0$ so
$T^{i-1}$ is $K_4$ with one edge removed, \ie\ the graph $(1)$ in Figure \ref{K4nei}.
If $T = T^i$ then $T^{i-1}$ is a connected subgraph, so using the first case we can replace $T$ by $T^{i-1}$ (thus avoiding $K_4$).
Otherwise $T^{i+1}$ must be the graph $(3)$ in Figure \ref{K4nei}.
Then we can replace
$T^i$ by the graph $(4)$ in that Figure, proving the claim.
\begin{figure}[!h]
\ifXY
\begin{equation}\label{avoidK4}\nonumber
\xymatrix@C=10pt@R=10pt{
\circ \ar@{-}[r] \ar@{-}[d]
    & \circ \ar@{-}[ld] \ar@{-}[d]
    & \circ \ar@{-}[r] \ar@{-}[rd]\ar@{-}[d]
    & \circ\ar@{-}[ld]   \ar@{-}[d]
    & \circ \ar@{-}[r] \ar@{-}[rd] \ar@{-}[d]
    & \circ \ar@{-}[r] \ar@{-}[ld]\ar@{-}[d]
    & \circ
    & \circ \ar@{-}[r] \ar@{-}[d]
    & \circ \ar@{-}[r] \ar@{-}[ld] \ar@{-}[d]
    & \circ \\
\circ \ar@{-}[r]_{(1)}
    & \circ
    & \circ \ar@{-}[r]_{(2)}
    & \circ
    & \circ \ar@{-}[r]_{(3)}
    & \circ
    &
    & \circ \ar@{-}[r]_{(4)}
    & \circ
    &       \\
}
\end{equation}
\else
\bigskip

(graphs around $K_4$).

\bigskip
\fi 
\caption{}\label{K4nei}
\end{figure}

We continue under the assumption that $T'$ is a disconnected subgraph of the connected graph $T$, with $|T-T'| = 1$.
In order to show that $\CoxY{T'} \sub \CoxY{T}$, we prove that
$F_{t,n-1}$ embeds in $F_{t,n}$. Consider the commutative diagram
\ifXY
\begin{equation}\nonumber
\xymatrix{
1 \ar[r]
    & F_{t,n} \ar@{^{(}->}[r]
    & \Fstar \ar[r]^{\ab}
    & \Z^{t} \ar@{=}[d] \ar[r]
    & 1 \\
1 \ar[r]
    & F_{t-1,n} \ar@{^{(}->}[r] \ar[u]^{\rho}
    & \Fstar[t-1,n] \ar[r]^{\ab} \ar[u]^{\rho}
    & \Z^t  \ar[r]
    & 1 \\
}
\end{equation}
\else

\bigskip

(A commutative diagram for $F_{t,n-1} \sub F_{t,n}$).

\bigskip
\fi 
where $\rho$ is defined by $x_i \mapsto x_i$ ($i = 1,\dots,n$).
This induces an embedding $F_{t-1,n} \hookrightarrow F_{t,n}$
which agrees with the $S_n$ action.

Letting $S_{n-1}$ denote the stabilizer in $S_n$ of the isolated vertex in $T'$, we have
$$\CoxY{T'} \isom \semidirect{S_{n-1}}{F_{t,n-1}} \hookrightarrow \semidirect{S_{n-1}}{F_{t,n}} \hookrightarrow \semidirect{S_{n}}{F_{t,n}} \isom \CoxY{T},$$
and the image of $\CoxY{T'}$ is easily computed to be the subgroup $\sg{T'}$ of $\CoxY{T'}$.
\end{proof}

The proof of Proposition \ref{At1At} shows that $A_{t-1,n}$ is a quotient group of $A_{t,n}$, so
we have proved
\begin{rem}\label{T1quotT}
If $T' \sub T$ is connected, then $\CoxY{T'}$ is a retract of $\CoxY{T}$.
\end{rem}
If $T = T'\cup \set{x}$, then the group $\CoxY{T'}$ is obtained by
adding to $\CoxY{T}$ the cyclic relation associated to the basic
cycle of $x$. In other words, assuming the cyclic relation of a
basic cycle amounts to erasing the non-$T_0$ edge on that cycle
from $T$.

Of course Remark \ref{T1quotT} does not hold if $T'$ is not
connected, as $S_{n_1} \times S_{n_2}$ is not a quotient of
$S_{n_1+n_2}$.

Until now we always assumed that $T$ is connected. The generators from distinct components commute, so we have
\begin{rem}
Let $T = T_1 \cup \dots \cup T_s$ be the decomposition to connected components of the graph $T$.

Then $\Cox{T} = \Cox{T_1} \times \dots \times \Cox{T_s}$,
$\CoxY{T} = \CoxY{T_1} \times \dots \times \CoxY{T_s}$, and $K(T)
= K(T_1) \times \dots \times K(T_s)$. The map $\maptoSn \co
\Cox{T} \ra S_n$ covers $S_{n_1} \times \dots \times S_{n_s}$
where $n_i$ is the number of vertices in $T_i$.
\end{rem}

\subsection{The Coxeter graph}\label{CoxGraph}

It is customary to define a Coxeter group by a graph (called the
Coxeter graph of the group), defined on the set of generators,
where generators $u,v$ are connected in the graph iff they do not
commute; if the order of $uv$ is assumed to be a number $p > 3$,
the edge connecting $v,u$ is labelled by $p$.

{}From our definition of $\Cox{T}$ it follows that the Coxeter
graph of $\Cox{T}$ is the dual graph $T^{\#}$, which is defined on
the set of edges of $T$, and two edges of $T$ are connected in
$T^{\#}$ iff they intersect in $T$.

\begin{exmpl}
Let $T$ be a cycle on $n$ vertices: since all the vertices are of degree $2$,
$\Cox{T} = \CoxY{T}$. The Coxeter graph is again a cycle on $n$ vertices,
being the dual graph of $T$.
We obtain the well known result
, that the Coxeter group of a cycle on $n$ vertices, is isomorphic
to $\semidirect{S_n}{\Z^{n-1}}$.
\end{exmpl}

\begin{exmpl}
Let $Y$ be the graph on four vertices, consisting of three edges
meeting in a vertex, and $\Delta$ the graph of a triangle. Then
the dual of both graphs is a triangle, showing that $\Cox{Y} =
\Cox{\Delta}$. Moreover, since in $\Delta$ only two edges meet in
every vertex, $\CoxY{\Delta} = \Cox{\Delta}$. So we have that
$\Cox{Y} = \Cox{\Delta} = \CoxY{\Delta} = \semidirect{S_3}{\Z^2}$,
while $\CoxY{Y} = S_4$. Let $u,v,w$ denote the edges of $Y$. Then
in the isomorphism $\Cox{Y} = \semidirect{S_3}{\Z^2}$ we can take
$u,v$ as generators of the $S_3$ piece, and $\sg{uvwv, vuwu}$ as
generators of $\Z^2$. The kernel of the map to $\CoxY{Y} = S_4$ is
$\sg{[u,vwv],[v,uwu]} \isom \Z^2$.
\end{exmpl}

The graphs $T^{\#}$ tend to have many edges. In particular
not every graph is of the form $T^{\#}$.
The easiest example is the following.
\begin{rem} 
If a graph $S$ has the graph of Figure \ref{forkgr} as a subgraph,
then it is not of the form $T^{\#}$ for any graph $T$.
\begin{figure}[!h]
\ifXY
\begin{equation}\label{forkdi}\nonumber
\xymatrix@C=10pt@R=10pt{
{}
    & \circ \ar@{-}[d]
    &
\\
\circ \ar@{-}[r]
    & \circ \ar@{-}[r]
    & \circ
\\
}
\end{equation}
\else

\bigskip

(A picture of a fork).

\bigskip
\fi 
\caption{A forbidden subgraph of $T^{\#}$}\label{forkgr}
\end{figure}
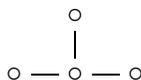
\end{rem}

Using this observation, one easily scans the lists of finite or
hyperbolic Coxeter groups (which can be found in \cite{Humph}, for
example), to check the following
\begin{rem}
The only hyperbolic Coxeter groups of the form $\Cox{T}$ are the symmetric groups $S_n$, and the three infinite
non-compact groups obtained for the graphs of Figure \ref{hyperT}.
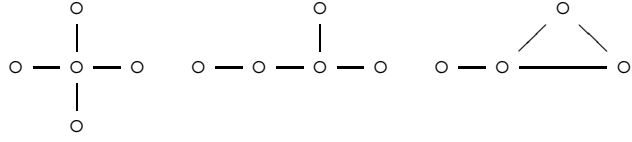
\begin{figure}[!h]
\ifXY
\begin{equation}\label{hyperbolicT}\nonumber
\xymatrix@C=10pt@R=10pt{
{}
    & \circ \ar@{-}[d]
    &
    &
    &
    & \circ \ar@{-}[d]
    &
    &
    &
    & \circ\ar@{-}[ld] \ar@{-}[rd]
    & \\
 \circ\ar@{-}[r]
    & \circ\ar@{-}[r]\ar@{-}[d]
    & \circ
    & \circ\ar@{-}[r]
    & \circ\ar@{-}[r]
    & \circ\ar@{-}[r]
    & \circ
    & \circ\ar@{-}[r]
    & \circ\ar@{-}[rr]
    &
    & \circ \\
{}
    & \circ
    &
    &
    &
    &
    &
    &
    &
    &
    & \\
}
\end{equation}
\else

\bigskip

(The three $T$ graphs).

\bigskip
\fi 
\caption{The hyperbolic $T$s}\label{hyperT}
\end{figure}

The corresponding (dual) Coxeter graphs are shown in Figure \ref{hyperTdu}.
The groups $\CoxY{T}$ are $S_5$, $S_5$ and $\semidirect{S_4}{\Z^3}$,
respectively.

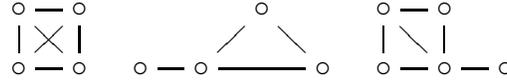
\begin{figure}[!h]
\ifXY
\begin{equation}\label{hyperbolicduals}\nonumber
\xymatrix@C=10pt@R=10pt{
\circ \ar@{-}[d] \ar@{-}[r] \ar@{-}[rd] & \circ \ar@{-}[d] \ar@{-}[ld]
    &
    &
    & \circ\ar@{-}[ld]\ar@{-}[rd]
    &
    & \circ\ar@{-}[r]\ar@{-}[d]\ar@{-}[rd]
    & \circ\ar@{-}[d]
    &  \\
\circ \ar@{-}[r]
    & \circ
    & \circ \ar@{-}[r]
    & \circ \ar@{-}[rr]
    &
    & \circ
    & \circ \ar@{-}[r]
    &\circ\ar@{-}[r]
    &\circ  \\
}
\end{equation}
\else

\bigskip

(The three $T^{\#}$ graphs).

\bigskip
\fi 
\caption{Hyperbolic $T^{\#}$s}\label{hyperTdu}
\end{figure}
\end{rem}

\subsection{Coxeter groups defined by signed graphs}\label{others}

Let $V$ be the vector space spanned by the generators of a Coxeter
group. In the Tits' representation of Coxeter groups, a generator
$u$ acts as $r_u \co V \ra V$, where $r_u$ is defined by $r_u(u) =
-u$, $r_u(v) = v$ if $u,v$ commute in the group, and $r_u(v) =
v+u$ if $(uv)^3 = 1$. Changing this slightly, to $r_u(v) = v-u$ if
$(uv)^3 = 1$, we obtain a representation of $C_Y(T)$.

This observation, in a more general context, is one of the
motivations of \cite{CST} to introduce certain quotients of
Coxeter groups, indexed by signed graphs. We direct the reader to
\cite{CST} for the definitions and notation used here. Let $f$ be
a signing of the graph $T^{\#}$. If for every $u,v,w \in T$ which
meet in a vertex, the resulting triangle in $T^{\#}$ is odd (with
respect to $f$), then $\operatorname{Cox}(T{^\#},f)$ is a quotient
of $\CoxY{T}$. It is fairly direct to evaluate the 'cut elements'
(which define $\operatorname{Cox}(T{^\#},f)$ as a quotient of
$\Cox{T}$) in $\semidirect{S_n}{A_{t,n}}$, so these groups can be
computed using our results. For example, the group $D_2$ computed
in \cite{FJNT} is our $\CoxY{T}$ for $T$ the graph of Figure
\ref{K4nei}.(1) (indeed
$D_2 = \semidirect{S_4}{F_{2,4}}$). 

Tsaranov groups are generated by elements $t_1,\dots,t_n$
of order $3$ with relations $(t_i t_j^{\pm 1})^3 = 1$, where the signs are
indexed by a graph
(the exponent being $-1$ iff $(i,j)$ is in the graph).
In \cite{CST} it is shown that
$\Tsar$ (which is the Tsaranov group of
the graph $\Gamma$, adjoint with an involution whose conjugation
inverts the generators) is a certain quotient of the Coxeter group
defined by the graph $\overline{\Gamma+0}$ (the complement of a
disjoint union of $\Gamma$ and a point $0$).

\begin{prop}
Let $a,b \geq t \geq 0$. Let $K_{a,b}$ be the complete bipartite
graph on $a+b$ vertices, and let $\Gamma$ be the graph obtained by deleting
$t$ disjoint edges from $K_{a,b}$.
Set $n = a+b+2-t$, and let
$X = \set{x,y,\dots}$ be a set of size $t$ indexing the generators of $F_{t,n}$.

Then $\Tsar$ is the subgroup $\semidirect{S_{n}}{F_{t,n}}$ of
$\semidirect{S_{n}}{\F_t^{n}}$, modulo its normal
subgroup $\sg{x_{i}^2x_{j}^{-2}}^{F_{t,n}}$.
\end{prop}
\begin{proof}
Taking $T$ to be the union of $t$ triangles with a common vertex
$u_0$, to which we glue a star of $a-t$ edges at one end and a star
of $b-t$ edges at the other, we see that $\overline{\Gamma+0} =
T^{\#}$.

>From Theorem 8.3 of \cite{CST} it then follows that $\Tsar$ is a quotient of
$\CoxY{T}$. The presentation given in Theorem 8.1 of that paper
translates to adding the relations $x_{ij}^2 = 1$ to $\CoxY{T} = \semidirect{S_n}{A_{t,n}}$. Since $T$ can never be the graph $K_4$, Proposition \ref{AintoAstar} applies, and $A_{t,n} \sub \Astar$, with $\sg{x_{ij}^2}^{A_{t,n}}$
being a normal subgroup in $\Astar$. Mapping this presentation
to $\Fstar$, we get the desired result.
\end{proof}

Of special interest is the Tsaranov group of a hexagon
(see \cite[p. 179]{FJNT}). Note that deleting
three disjoint edges from $K_{3,3}$ gives a hexagon, so we get
\begin{cor}\label{TsaranovHexagon}
The Tsaranov group of a hexagon is
$\semidirect{S_{5}}{F_{3,5}}$ modulo the relations
$x_{i}^2x_{j}^{-2}$ ($x \in X$).
\end{cor}

\end{document}